\documentclass[12pt]{article}
\usepackage[latin1]{inputenc} 
\usepackage{amsmath}
\usepackage{amssymb}
\usepackage{epsfig}
\usepackage{multirow}
\begin{document}
\begin{flushright}
	{ to appear in \emph{Discrete and Computational Geometry}}
\end{flushright}

\bigskip

\begin{center}

{\Large \bf Hamiltonian submanifolds of regular polytopes}

\medskip

\bigskip

{\sc Felix Effenberger} and {\sc Wolfgang K\"uhnel}
\end{center}

\bigskip

\bigskip
{\small {\sc Abstract:} 
We investigate polyhedral $2k$-manifolds as subcomplexes
of the boundary complex of a regular polytope. 
We call such a subcomplex {\it $k$-Hamiltonian}
if it contains the full $k$-skeleton of the polytope. 
Since the case of the cube is well known and since
the case of a simplex was also previously studied (these are 
so-called {\it super-neighborly triangulations}) 
we focus on the case of the
cross polytope and the sporadic regular 4-polytopes.
By our results the existence of 1-Hamiltonian surfaces 
is now decided for all regular polytopes.
Furthermore we investigate 2-Hamiltonian 4-manifolds
in the $d$-dimensional cross polytope. 
These are the ``regular cases'' satisfying equality in Sparla's inequality.
In particular, we present a new example with 16 vertices which is
highly symmetric with an automorphism group of order 128. 
Topologically it is homeomorphic to a connected sum of 7 copies
of $S^2 \times S^2$.
By this example all regular cases of $n$ vertices 
with $n < 20$ or, equivalently, all cases of regular
$d$-polytopes with $d\leq 9$ are now decided.


} 

\subsection*{1. Introduction and results}
The idea of a Hamiltonian circuit in a graph can be generalized to
higher-dimensional complexes as follows: A subcomplex $A$ of a
polyhedral complex $K$ is called {\sf $k$-Hamiltonian}\footnote{not to
be confused with the notion of a $k$-Hamiltonian 
graph \cite{KT}} if $A$ contains
the full $k$-dimensional skeleton of $K$. 
It seems that this concept was first developed by 
C.Schulz \cite{Sch1, Sch2}. A Hamiltonian circuit
then becomes a special case of a 0-Hamiltonian subcomplex 
of a 1-dimensional graph or of a higher-dimensional complex \cite{E}.
If $K$ is the boundary complex of a convex polytope then this
concept becomes particularly interesting and quite geometrical
\cite[Ch.3]{Ku1}.
A.Altshuler \cite{A} investigated 1-Hamiltonian closed surfaces 
in special polytopes.
A triangulated surface with a complete edge graph $K_n$
can be regarded as a 1-Hamiltonian subcomplex of the simplex with
$n$ vertices. These are the so-called {\it regular cases} in
Heawood's Map Color Theorem \cite{R}, \cite[2C]{Ku1}, and people talk
about the uniquely determined {\it genus of the complete graph $K_n$} 
which is (in the orientable regular cases $n\equiv 0,3,4,7\ (12),\, n\geq 4$)
$$g = \frac{1}{6}{{n-3} \choose 2}.$$
Moreover, the induced piecewise linear embedding of the surface
into Euclidean $(n-1)$-space then has the two-piece property,
and it is tight \cite[2D]{Ku1}.

\medskip
Centrally-symmetric analogues can be regarded as 1-Hamiltonian
subcomplexes of cross polytopes or other centrally symmetric polytopes, 
see \cite{Ku2}.
Similarly we have the {\it genus of the $d$-dimensional
cross polytope} \cite{JR} which is (in the orientable regular cases 
$d\equiv 0,1\ (3),\, d\geq 3$)
$$g = \frac{1}{3}(d-1)(d-3).$$
There are famous examples of quadrangulations of surfaces 
originally due to H.~S.~M.~Coxeter
which can be regarded as 1-Hamiltonian subcomplexes of higher-dimensional cubes
\cite {KS}, \cite[2.12]{Ku1}. 
Accordingly one talks about the {\it genus of the $d$-cube} 
(or rather its edge graph) which is (in the orientable case)
$$g = 2^{d-3}(d-4) + 1,$$ see \cite{R0}, \cite{BH}.
However, in general the genus of a 1-Hamiltonian surface in
a convex $d$-polytope is not uniquely determined, as pointed out
in \cite{Sch1, Sch2}. This uniqueness
seems to hold especially for regular polytopes where the regularity
allows a computation of the genus by a simple counting argument.

\medskip
In the cubical case there are higher-dimensional
generalizations by Danzer's construction
of a {\it power complex} $2^K$ for a given simplicial complex $K$.
In particular there are many examples of $k$-Hamiltonian $2k$-manifolds
as subcomplexes of higher-dimensional cubes, see \cite{KS}.
For obtaining them one just has to start with a neighborly simplicial
$(2k-1)$-sphere $K$.
A large number of the associated complexes $2^K$ are topologically connected 
sums of copies of $S^k \times S^k$. This seems to be the standard case.

\medskip
Concerning triangulations of manifolds, a $d$-dimensional simplicial complex
is called a {\sf combinatorial $d$-manifold} if the union of its simplices
is homeomorphic to a $d$-manifold and if the link of each $k$-simplex
is a combinatorial $(d-k-1)$-sphere. In what follows all
triangulations of manifolds are assumed to be combinatorial.
There exist triangulations 
of manifolds which are not combinatorial, for an example
based on the Edwards sphere see \cite{Bj-Lu}.

\medskip
With respect to the simplex as the ambient polytope a $k$-Hamiltonian
subcomplex is also called a {\sf $(k+1)$-neighborly triangulation} since
any $k+1$ vertices are common neighbors in a $k$-dimensional
simplex. The crucial case is the case of $(k+1)$-neighborly
triangulations of $2k$-manifolds. This case was studied by the second
author in \cite{Ku1}. One could call this the case
of {\sf super-neighborly triangulations} in analogy with neighborly polytopes:
The boundary complex of a $(2k+1)$-polytope can be at most 
$k$-neighborly unless it is a simplex.
However, combinatorial $2k$-manifolds can go beyond
$k$-neighborliness, depending on the topology.
Except for the trivial case of the boundary
of a simplex itself there are only a finite number of known examples
of super-neighborly triangulations,
reviewed in \cite{K-Lu}. They are necessarily tight \cite[Ch.4]{Ku1}, compare 
Section 5 below.
The most significant ones are
the unique 9-vertex triangulation of the complex projective plane \cite{Ku-Ba},
\cite{Ku-La1}, a 16-vertex triangulation of a K3 surface \cite{CK}
and several 15-vertex triangulations of an 8-manifold ``like the quaternionic
projective plane'' \cite{Br-Ku2}. There is also an asymmetric
13-vertex triangulation of $S^3 \times S^3$, but most of the examples
are highly symmetric.
For any $n$-vertex triangulation of a $2k$-manifold $M$
the generalized Heawood inequality
$${{n-k-2} \choose {k+1}} \geq {{2k+1} \choose {k+1}}(-1)^k(\chi(M) - 2)$$
was conjectured in \cite{Ku0}, \cite{Ku1} and later almost completely proved
by I.~Novik in \cite{N1} and proved in \cite{NS}.
Equality holds precisely in the case of super-neighborly triangulations.
These are $k$-Hamiltonian in the $(n-1)$-dimensional simplex.
In the case of 4-manifolds (i.e., $k=2$) an elementary proof was
already contained in \cite[4B]{Ku1}.

\medskip
In the case of 2-Hamiltonian subcomplexes of cross polytopes the
first non-trivial example was constructed by E.~Sparla as a centrally-symmetric
12-vertex triangulation of $S^2 \times S^2$ as a subcomplex
of the boundary of the 6-dimensional cross polytope \cite{Sp1}, \cite{La-Sp}.
Sparla also proved the following analogous Heawood inequality for
the case of 2-Hamiltonian 4-manifolds in centrally symmetric $d$-polytopes
$${{\frac{1}{2}(d-1)} \choose {3}} \leq 10(\chi(M) - 2)$$
and the opposite inequality for centrally-symmetric triangulations
with $n = 2d$ vertices.

\medskip
Higher-dimensional examples were found by F.~H.~Lutz \cite{Lu}:
There are centrally-symmetric 16-vertex triangulations of $S^3 \times S^3$ 
and 20-vertex triangulations of $S^4 \times S^4$.
The 2-dimensional example in this series is the well known 
unique centrally-symmetric 8-vertex torus \cite[3.1]{Ku2}.
All these are tightly embedded into the ambient Euclidean space \cite{K-Lu}.
The generalized Heawood inequality for centrally symmetric $2d$-vertex
triangulations of $2k$-manifolds 
$$4^{k+1}{{\frac{1}{2}(d-1)} \choose {k+1}} \geq {{2k+1} \choose {k+1}}(-1)^k(\chi(M) - 2)$$
was conjectured
by Sparla in \cite{Sp2} and later almost completely proved by I.~Novik
in \cite{N2}. 

\bigskip
In the present paper we show that Sparla's inequality for 2-Hamiltonian
4-manifolds in the skeletons of $d$-dimensional cross polytopes
is sharp for $d \leq 9$. More precisely, we show that each of the
regular cases (that is, the cases of equality) for $d \leq 9$ really occurs. 
Since the cases $d=7$
and $d=9$ are not regular, the crucial point is
the existence of an example for $d=8$ and, necessarily, $\chi = 16$.  
In addition we examine the case of $1$-Hamiltonian surfaces
in the three sporadic regular 4-polytopes, see Section 2.
It seems that so far no decision about existence or non-existence
could be made, compare \cite{S}.

\subparagraph{Main Theorem}
\begin{enumerate}
\item {\it 
All cases of $1$-Hamiltonian surfaces in the regular polytopes
are decided.
In particular there are no $1$-Hamiltonian surfaces in the $24$-cell,
$120$-cell or $600$-cell.} 
\item {\it
All cases of $2$-Hamiltonian
$4$-manifolds in the regular $d$-polytopes are decided up to
dimension $d=9$.
In particular, there is a new example of a
$2$-Hamiltonian
$4$-manifold in the boundary complex of the $8$-dimensional cross polytope.}
\end{enumerate}

\medskip
This follows from certain known results and a 
combination of Propositions 1, 2, 3, and Theorem 2 below.

\medskip
The regular cases of $1$-Hamiltonian surfaces are the following,
and each case occurs:

\qquad
\begin{tabular}{lll}
$d$-simplex: & $d\equiv 0,2 \ (3)$ &\cite{R}\\
$d$-cube: &any $d\geq 3$&\cite{BH},\cite{R0}\\
$d$-octahedron: &$d\equiv 0,1 \ (3)$& \cite{JR}.\\
\end{tabular}

\medskip
The regular cases of $2$-Hamiltonian $4$-manifolds for $d\leq 9$
are the following:

\qquad
\begin{tabular}{lll}
$d$-simplex:  & $d=5,8,9$&\cite{Ku-La1}\\
$d$-cube: &$d=5,6,7,8,9$&\cite{KS}\\
$d$-octahedron: &$d=5,6,8$&Theorem 2.
\end{tabular}

Here each of these cases occurs except for the case of the $9$-simplex
\cite{Ku-La1}. 
Furthermore $2$-Hamiltonian $4$-manifolds in the $d$-cube are known to exist
for any $d\geq 5$ \cite{KS}. In the case of the $d$-simplex the 
next regular case $d=13$
is undecided, the case $d=15$ occurs \cite{CK}.
The next regular case of a $d$-octahedron is the case $d=10$, see Remark 2
below.
%%%%%%%%%%%%%%%%%%%%%%%% SECTION 2  %%%%%%%%%%%%%%%%%%%%%%%%%%%%%%%%%%%%%%%
\subsection*{2. Hamiltonian surfaces in the $24$-cell, $120$-cell, $600$-cell}
There are Hamiltonian cycles in each of the Platonic solids.
The numbers of distinct Hamiltonian cycles (modulo symmetries of the
solid itself) are $1,1,2,1,17$ for the cases of the tetrahedron, cube,
octahedron, dodecahedron, icosahedron, see \cite[pp.\ 277 ff.]{H}. 
A 1-Hamiltonian surface in the boundary complex of a Platonic
solid must coincide with the boundary itself and is, therefore, not 
really interesting.

\medskip
Hamiltonian cycles in the regular 4-polytopes are known to exist. However,
it seems that 1-Hamiltonian surfaces in the 2-skeleton of any of the
three sporadic regular 4-polytopes have not yet been systematically
investigated.
A partial attempt can be found in \cite{S}.
\subsubsection*{2.1 The $24$-cell}
The boundary complex of the $24$-cell $\{3,4,3\}$ consists of 24
vertices, 96 edges, 96 triangles and 24 octahedra.
Any 1-Hamiltonian surface (or pinched surface) 
must have 24 vertices, 96 edges and, consequently,
64 triangles, hence it has Euler characteristic $\chi = -8$.
Every edge in the polytope is in three triangles.
Hence we must omit exactly one of them in each case for getting a surface where
every edge is in two triangles. 
Since the vertex figure in the polytope is a cube, 
each vertex figure in the surface
is a Hamiltonian circuit of length 8 in the edge graph of a cube.
It is well known that this circuit is uniquely determined up to
symmetries of the cube.
Starting with one such vertex figure, there are four missing edges in
the cube which, therefore, must be in the uniquely determined other 
triangles of the 24-cell. In this way, one can inductively construct an
example or, alternatively, verify the non-existence.
If singular vertices are allowed, then the only possibility
is a link which consists of two circuits of length four each.
This leads to the following proposition.

\begin{table}[hbt]
\centering
{\small
\begin{tabular}{@{}|c|c|c|c|@{}} \hline
type&group&order&generators\\\hline

\multirow{2}{*}{1}&
\multirow{2}{*}{$C_4 \times C_2$}&
\multirow{2}{*}{8}&

{\scriptsize $( 1\,12\,16\,18)( 2\,17\,23\, 7)( 3\,13\,20\,21)( 4\,22\,11\, 5)( 6\,19)( 8\,24\,14\,10)$,}\\
&&&{\scriptsize $( 1\, 3)( 4\, 8)( 5\,10)( 9\,15)(11\,14)(12\,13)(16\,20)(18\,21)(22\,24)$}\\\hline

\multirow{2}{*}{2}&
\multirow{2}{*}{$D_8$}&
\multirow{2}{*}{8}&
{\scriptsize $( 1\,16)( 2\,17)( 3\,22)( 5\,20)( 6\, 9)( 7\,23)( 8\,12)(10\,24)(14\,18)(15\,19)$,}\\
&&&{\scriptsize $( 2\, 3)( 4\, 6)( 5\, 7)( 9\,11)(12\,14)(13\,15)(17\,20)(19\,21)(22\,23)$}\\\hline

\multirow{2}{*}{3}&
\multirow{2}{*}{$C_2\times C_2$}&
\multirow{2}{*}{4}&
{\scriptsize $( 1\,24)( 2\,13)( 3\,15)( 4\,17)( 5\,19)( 6\,20)( 7\,21)( 9\,22)(11\,23)$,}\\
&&&{\scriptsize $( 2\, 5)( 3\, 7)( 4\, 9)( 6\,11)( 8\,18)(13\,19)(15\,21)(17\,22)(20\,23)$}\\\hline

\multirow{2}{*}{4}&
\multirow{2}{*}{{\scriptsize $(((C_4 \times C_2):C_2):C_2):C_2$}}&
\multirow{2}{*}{64}&
{\scriptsize $( 1\, 8\,10\,12)( 3\,13\, 5\, 4)( 6\,15\,19\, 9)( 7\,17)(11\,20\,21\,22)(14\,24\,18\,16)$,}\\
&&&{\scriptsize $( 2\, 3)( 4\, 6)( 5\, 7)( 9\,11)(12\,14)(13\,15)(17\,20)(19\,21)(22\,23)$}\\\hline

\multirow{2}{*}{5}&
\multirow{2}{*}{$S_3$}&
\multirow{2}{*}{6}&
{\scriptsize $( 1\, 3)( 4\, 8)( 5\,10)( 9\,15)(11\,14)(12\,13)(16\,20)(18\,21)(22\,24)$,}\\
&&&{\scriptsize $( 1\,22\,15)( 2\,12\,13)( 3\, 9\,24)( 4\,17\, 8)( 5\,19\,10)( 6\,16\,20)( 7\,18\,21)(11\,23\,14)$}\\\hline

\multirow{3}{*}{6}&
\multirow{3}{*}{$C_2\times D_8$}&
\multirow{3}{*}{16}&
{\scriptsize $( 1\,11)( 2\,23)( 3\,14)( 4\,16)( 5\,18)( 8\,20)(10\,21)(12\,22)(13\,24)$,}\\
&&&{\scriptsize $( 1\, 5)( 3\,12)( 4\,10)( 6\,19)( 7\, 9)( 8\,13)(11\,18)(14\,22)(15\,17)(16\,21)(20\,24)$,}\\
&&&{\scriptsize $( 1\, 3)( 4\, 8)( 5\,10)( 9\,15)(11\,14)(12\,13)(16\,20)(18\,21)(22\,24)$}\\\hline

\end{tabular}
}
\caption{Automorphism groups of the Hamiltonian
pinched surfaces in the $24$-cell}
\end{table}

\subparagraph{Proposition 1}
{\it There is no $1$-Hamiltonian surface in the $2$-skeleton of the
  $24$-cell. However, there are six combinatorial types of
strongly connected
$1$-Hamiltonian pinched surfaces with a number of pinch points
ranging between $4$ and $10$ and with the genus ranging 
between $g=3$ and $g=0$.
The case of the highest genus is a surface of genus three 
with four pinch points.
The link of each of the pinch points in any of these types
is the union of two circuits of length four.}

\bigskip
The six types and their automorphism groups are listed in Tables 1 and
2 where the labeling of the vertices of the 24-cell coincides
with the standard one in {\tt polymake} \cite{G-J}.

\bigskip

\begin{table}
\centering
{\small
\begin{tabular}{@{}|c|c|c|c|@{}} \hline
type&\# \ p.\ pts.&$g$&orbits\\\hline

1&%
10&%
0&%
{\scriptsize $\langle 1\,2\,3\rangle_4$, $\langle 1\,2\,4\rangle_8$, $\langle 1\,3\,6\rangle_4$, $\langle 1\,4\,9\rangle_8$, $\langle 1\,5\,7\rangle_8$, $\langle 1\,5\,9\rangle_8$, $\langle 1\,6\,11\rangle_8$, $\langle 1\,7\,11\rangle_8$, $\langle 2\,5\,10\rangle_4$, $\langle 4\,6\,8\rangle_4$}\\\hline

\multirow{2}{*}{2}&%
\multirow{2}{*}{10}&%
\multirow{2}{*}{0}&%
{\scriptsize $\langle 1\,2\,3\rangle_4$, $\langle 1\,2\,4\rangle_8$, $\langle 1\,4\,9\rangle_4$, $\langle 2\,3\,8\rangle_4$, $\langle 2\,4\,12\rangle_8$, $\langle 2\,5\,10\rangle_4$, $\langle 2\,5\,12\rangle_4$, $\langle 2\,8\,13\rangle_8$,}\\
&&&{\scriptsize $\langle 2\,10\,13\rangle_8$, $\langle 4\,6\,8\rangle_4$, $\langle 8\,13\,15\rangle_4$, $\langle 10\,13\,19\rangle_4$}\\\hline

\multirow{3}{*}{3}&%
\multirow{3}{*}{8}&%
\multirow{3}{*}{1}&%
{\scriptsize $\langle 1\,2\,3\rangle_4$, $\langle 1\,2\,4\rangle_4$, $\langle 1\,3\,6\rangle_4$, $\langle 1\,4\,9\rangle_2$, $\langle 1\,6\,11\rangle_2$, $\langle 2\,3\,8\rangle_4$, $\langle 2\,4\,12\rangle_4$, $\langle 2\,5\,10\rangle_2$, $\langle 2\,5\,12\rangle_2$,}\\
&&&{\scriptsize $\langle 2\,8\,13\rangle_2$, $\langle 2\,10\,13\rangle_2$, $\langle 3\,6\,14\rangle_4$, $\langle 3\,7\,10\rangle_2$, $\langle 3\,7\,14\rangle_2$, $\langle 3\,8\,15\rangle_2$, $\langle 3\,10\,15\rangle_2$, $\langle 4\,6\,8\rangle_4$,}\\
&&&{\scriptsize $\langle 4\,6\,16\rangle_4$, $\langle 4\,8\,17\rangle_2$, $\langle 4\,9\,16\rangle_2$, $\langle 4\,12\,17\rangle_2$, $\langle 6\,8\,20\rangle_2$, $\langle 6\,11\,14\rangle_2$, $\langle 6\,16\,20\rangle_2$}\\\hline

4&%
8&%
%&%
1&%
{\scriptsize $\langle 1\,2\,3\rangle_{32}$, $\langle 1\,2\,4\rangle_{32}$}\\\hline

\multirow{3}{*}{5}&%
\multirow{3}{*}{6}&%
\multirow{3}{*}{2}&%
{\scriptsize $\langle 1\,2\,3\rangle_3$, $\langle 1\,2\,4\rangle_6$, $\langle 1\,3\,6\rangle_3$, $\langle 1\,4\,9\rangle_3$, $\langle 1\,5\,7\rangle_6$, $\langle 1\,5\,9\rangle_3$, $\langle 1\,6\,11\rangle_6$, $\langle 1\,7\,11\rangle_6$, $\langle 2\,4\,12\rangle_3$,}\\
&&&{\scriptsize $\langle 2\,5\,10\rangle_3$, $\langle 2\,5\,12\rangle_3$, $\langle 4\,6\,8\rangle_3$, $\langle 4\,6\,16\rangle_3$, $\langle 4\,8\,17\rangle_1$, $\langle 5\,7\,18\rangle_3$, $\langle 5\,10\,19\rangle_1$, $\langle 6\,11\,16\rangle_3$,}\\
&&&{\scriptsize $\langle 7\,11\,14\rangle_3$, $\langle 7\,18\,21\rangle_1$, $\langle 11\,14\,23\rangle_1$}\\\hline

6&%
4&%
%&%
3&%
{\scriptsize $\langle 1\,2\,3\rangle_8$, $\langle 1\,2\,4\rangle_8$, $\langle 1\,3\,6\rangle_8$, $\langle 1\,4\,9\rangle_8$, $\langle 1\,5\,7\rangle_{16}$, $\langle 1\,6\,11\rangle_8$, $\langle 1\,7\,11\rangle_8$}\\\hline

\end{tabular}
}
\caption{Generating orbits of the 6 types of Hamiltonian
pinched surfaces in the $24$-cell}
\end{table}

Type 1 is a pinched sphere which is based on a subdivision of the boundary
of the rhombidodecahedron, see Figure 1 (left).
Type 4 is just a $(4 \times 4)$-grid square torus where each 
square is subdivided by an extra vertex, see Figure 1 (right).
These 16 extra vertices are identified in pairs, leading to the 8 pinch points.

\setcounter{figure}{1}
\begin{center}
\epsfig{file=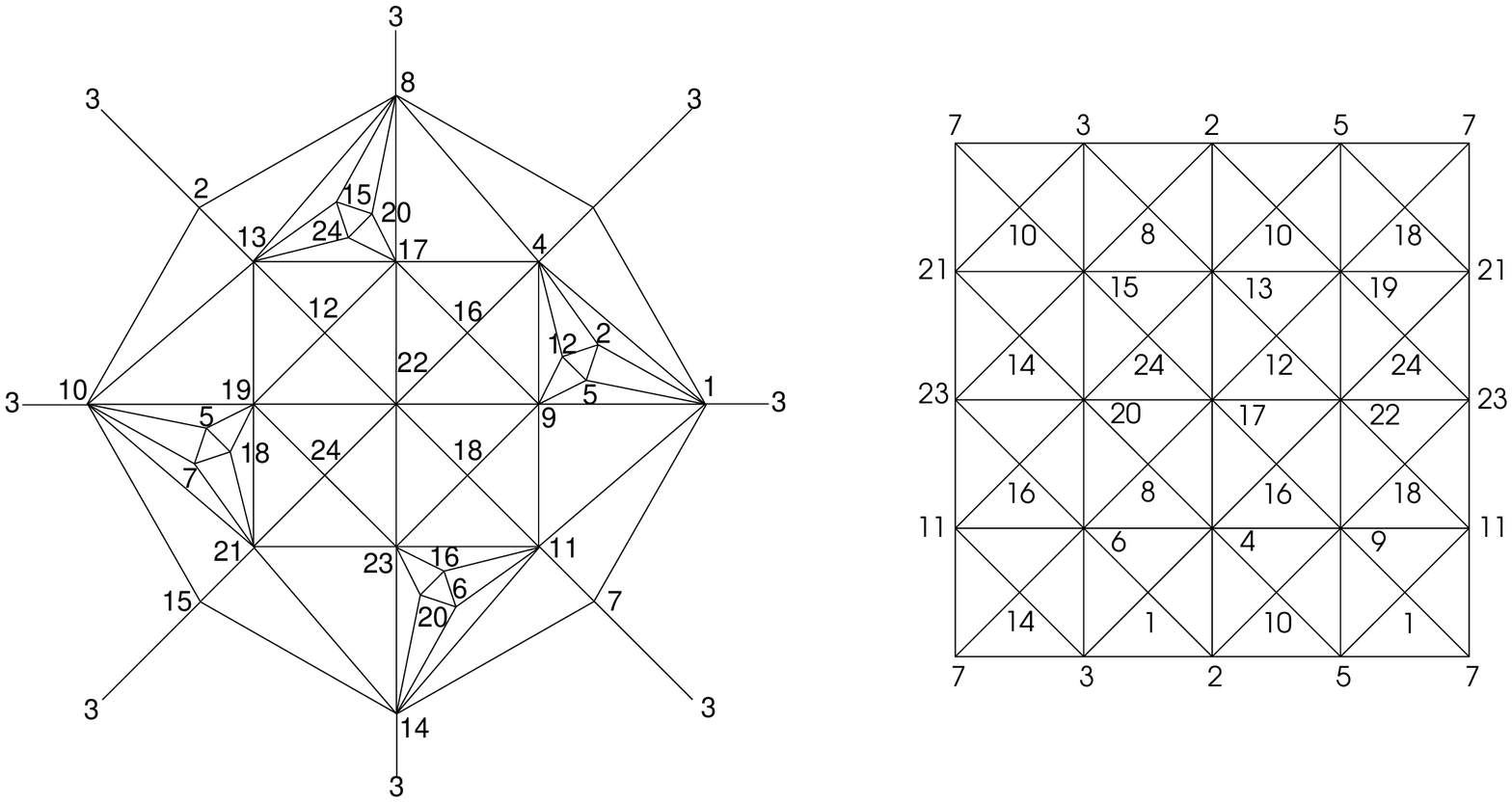,height=8.2cm}
Figure 1: Type 1 (left) and Type 4 (right) of
Hamiltonian pinched surfaces in the $24$-cell
\end{center}

\bigskip
Because $-8$ equals the Euler characteristic of the original 
(connected) surface
minus the number of pinch points it is clear that we can have at most
10 pinch points unless the surface splits into several components.
We present here in more detail
Type 6 as a surface of genus three with four pinch points, see Figure
3 (produced with JavaView). 
Its combinatorial type is given by the following list of 64 triangles:

\medskip
{\scriptsize
\centering

\begin{tabular}{@{}l@{}l@{}l@{}l@{}l@{}l@{}l@{}l@{}}
{$\langle 1\,2\,3 \rangle$},& {$\langle 1\,2\,4 \rangle$},& {$\langle 1\,3\,6 \rangle$},& {$\langle 1\,4\,9 \rangle$},& {$\langle 1\,5\,7 \rangle$},& {$\langle 1\,5\,9 \rangle$},& {$\langle 1\,6\,11 \rangle$},& {$\langle 1\,7\,11 \rangle$},\\
{$\langle 2\,3\,8 \rangle$},& {$\langle 2\,4\,8 \rangle$},& {$\langle 2\,5\,10 \rangle$},& {$\langle 2\,5\,12 \rangle$},& {$\langle 2\,10\,13 \rangle$},& {$\langle 2\,12\,13 \rangle$},& {$\langle 3\,6\,14 \rangle$},& {$\langle 3\,7\,10 \rangle$},\\
{$\langle 3\,7\,14 \rangle$},& {$\langle 3\,8\,15 \rangle$},& {$\langle 3\,10\,15 \rangle$},& {$\langle 4\,6\,8 \rangle$},& {$\langle 4\,6\,16 \rangle$},& {$\langle 4\,9\,12 \rangle$},& {$\langle 4\,12\,17 \rangle$},& {$\langle 4\,16\,17 \rangle$},\\
{$\langle 5\,7\,10 \rangle$},& {$\langle 5\,9\,18 \rangle$},& {$\langle 5\,12\,19 \rangle$},& {$\langle 5\,18\,19 \rangle$},& {$\langle 6\,8\,20 \rangle$},& {$\langle 6\,11\,14 \rangle$},& {$\langle 6\,16\,20 \rangle$},& {$\langle 7\,11\,18 \rangle$},\\
{$\langle 7\,14\,21 \rangle$},& {$\langle 7\,18\,21 \rangle$},& {$\langle 8\,13\,15 \rangle$},& {$\langle 8\,13\,17 \rangle$},& {$\langle 8\,17\,20 \rangle$},& {$\langle 9\,11\,16 \rangle$},& {$\langle 9\,11\,18 \rangle$},& {$\langle 9\,12\,22 \rangle$},\\
{$\langle 9\,16\,22 \rangle$},& {$\langle 10\,13\,19 \rangle$},& {$\langle 10\,15\,21 \rangle$},& {$\langle 10\,19\,21 \rangle$},& {$\langle 11\,14\,23 \rangle$},& {$\langle 11\,16\,23 \rangle$},& {$\langle 12\,13\,17 \rangle$},& {$\langle 12\,19\,22 \rangle$},\\
{$\langle 13\,15\,24 \rangle$},& {$\langle 13\,19\,24 \rangle$},& {$\langle 14\,15\,20 \rangle$},& {$\langle 14\,15\,21 \rangle$},& {$\langle 14\,20\,23 \rangle$},& {$\langle 15\,20\,24 \rangle$},& {$\langle 16\,17\,22 \rangle$},& {$\langle 16\,20\,23 \rangle$},\\
{$\langle 17\,20\,24 \rangle$},& {$\langle 17\,22\,24 \rangle$},& {$\langle 18\,19\,22 \rangle$},& {$\langle 18\,21\,23 \rangle$},& {$\langle 18\,22\,23 \rangle$},& {$\langle 19\,21\,24 \rangle$},& {$\langle 21\,23\,24 \rangle$},& {$\langle 22\,23\,24 \rangle$}. 
\end{tabular}

}

\bigskip
The pinch points are the vertices 2, 6, 19, 23 
with the following links:
{\small
\[
\begin{tabular}{rcc}
2:& (1 3 8 4) & (5 10 13 12)\\
6:& (1 3 14 11) & (4 8 20 16)\\
19:&(5 12 22 18) & (10 13 24 21)\\
23:&(11 14 20 16) & (18 21 24 22)
\end{tabular} 
\]}

The four vertices 7, 9, 15, 17  
are not joined to one another and
not to any of the pinch points either.
Therefore the eight vertex stars of $7,9,15,17,2,6,19,23$  
cover the 64 triangles of the surface entirely and simply,
compare Figure 2 where the combinatorial type is sketched.
In this drawing all vertices are 8-valent except for the four pinch points
in the two ``ladders'' on the right hand side which have to be
identified in pairs.

\medskip
The combinatorial automorphism group of order 16 is generated by
$$Z=(1\,11)(2\,23)(3\,14)(4\,16)(5\,18)(8\,20)(10\,21)(12\,22)(13\,24),$$
$$A=(1\,5)(3\,12)(4\,10)(6\,19)(7\,9)(8\,13)(11\,18)(14\,22)(15\,17)(16\,21)(20\,24),$$
$$B=(1\,3)(4\,8)(5\,10)(9\,15)(11\,14)(12\,13)(16\,20)(18\,21)(22\,24).$$

The elements $A$ and $B$ generate the dihedral group $D_8$ of order 8
whereas $Z$ commutes
with $A$ and $B$. Therefore the group is isomorphic with $D_8 \times C_2$.

\begin{figure}[hbt]
\centering
\epsfig{file=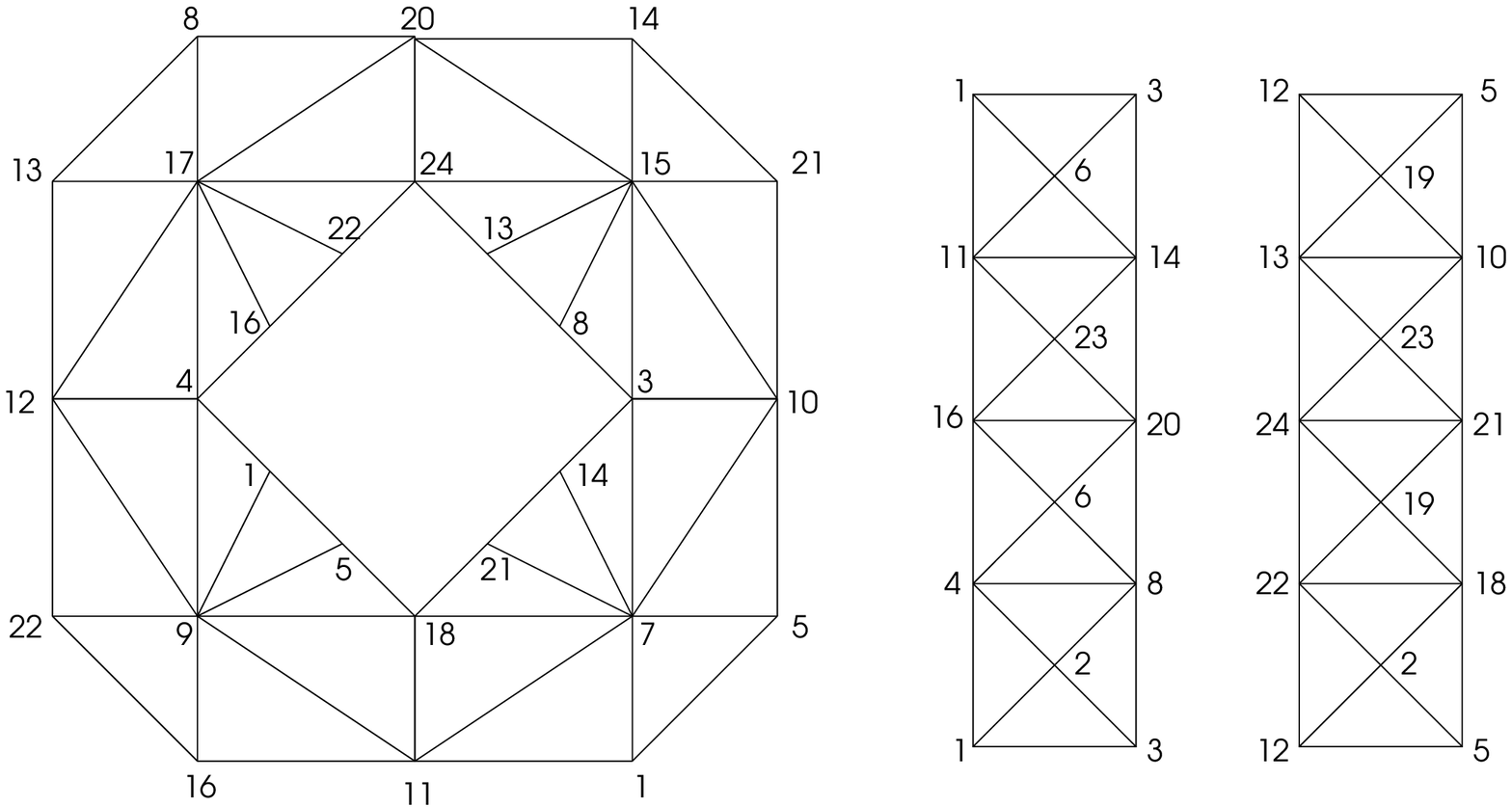,height=7.5cm}
\caption{The triangulation of the
Hamiltonian pinched surface of genus 3 in the $24$-cell}
\end{figure}

\begin{figure}[hbt]
\centering
\begin{tabular}{lr}
\epsfig{file=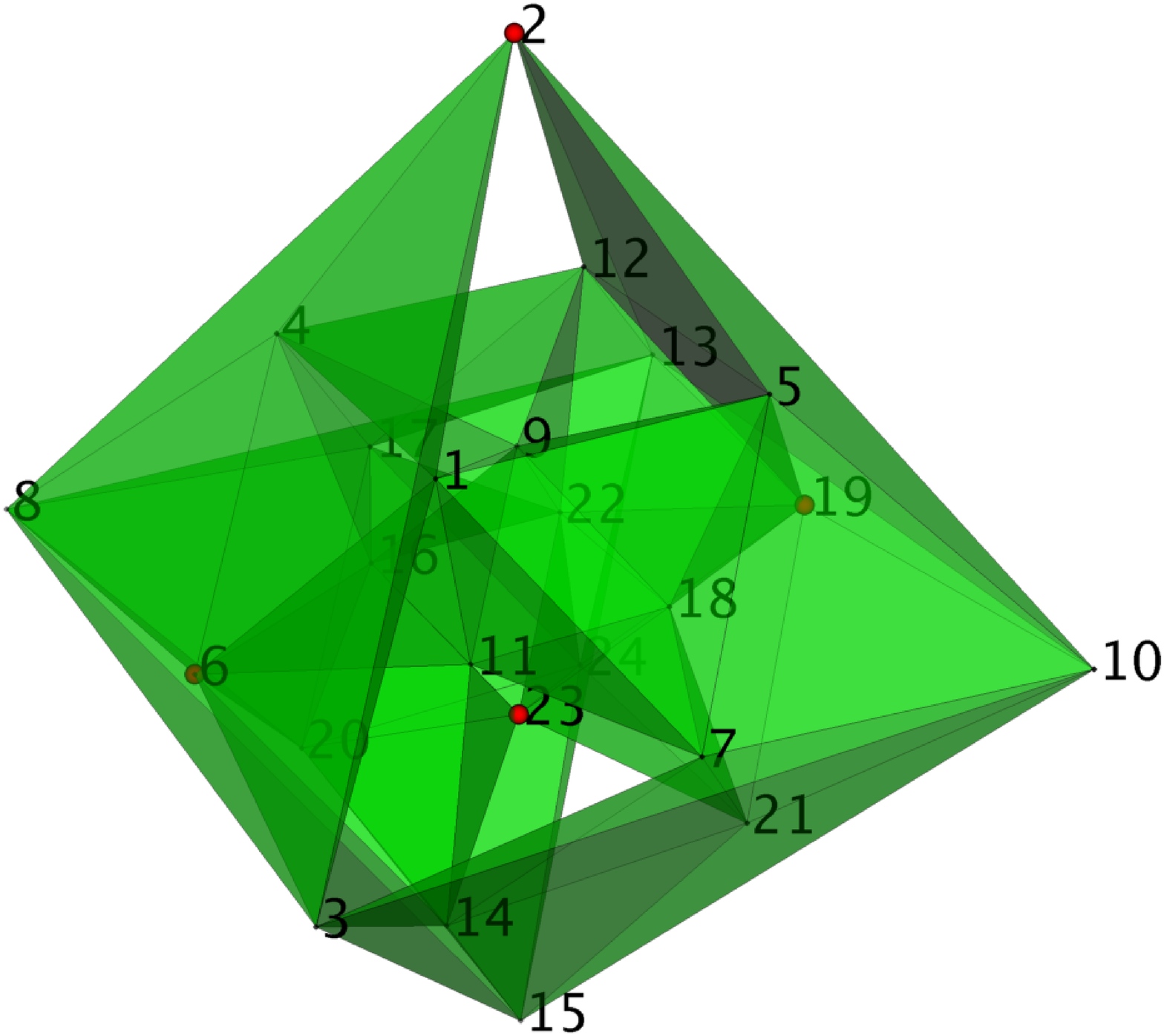,height=6.5cm}
&\epsfig{file=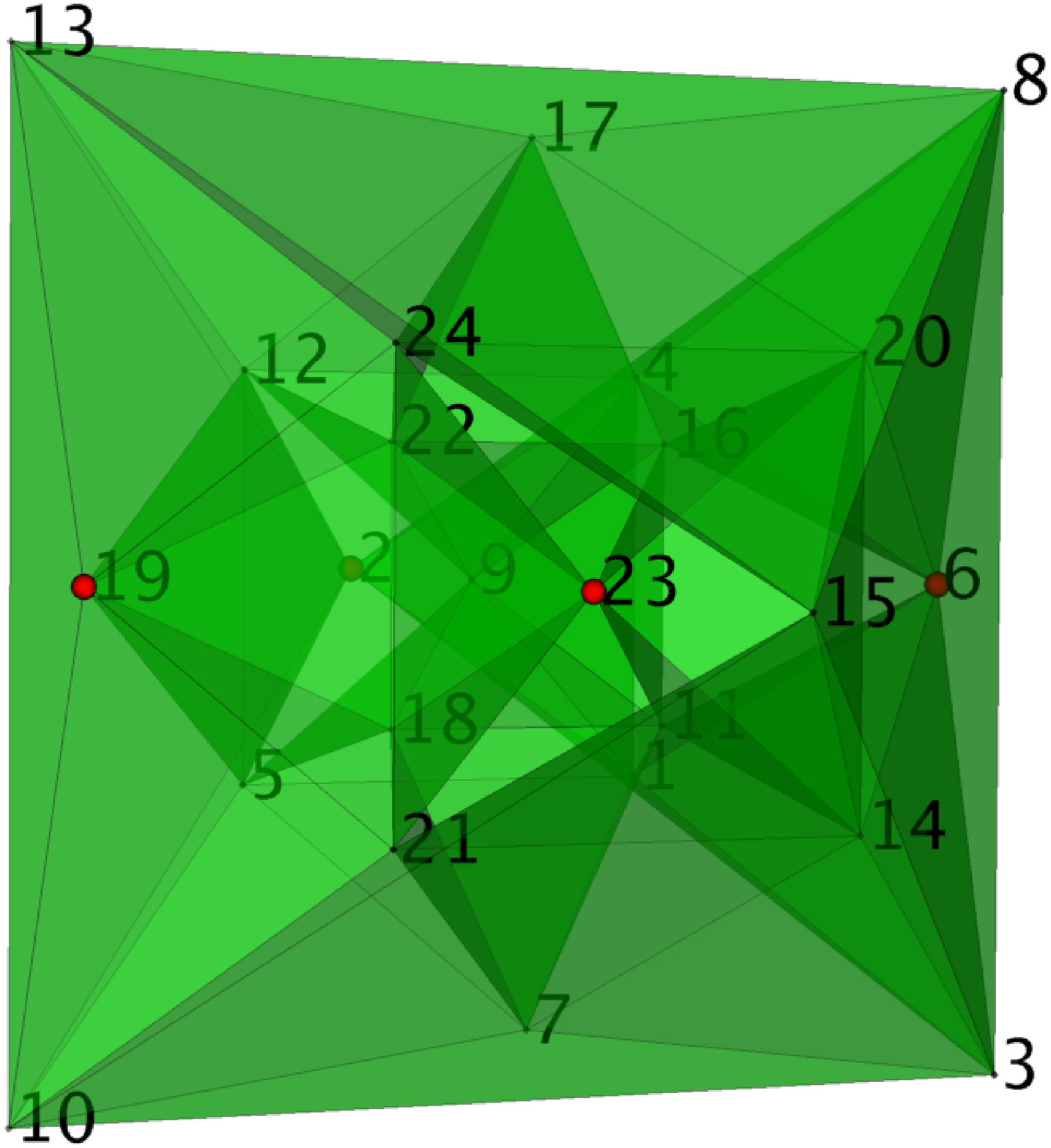,height=7cm}
\end{tabular}
\caption{Two projections of the
Hamiltonian pinched surface of genus 3 in the $24$-cell}
\end{figure}
\subsubsection*{2.2 The $120$-cell and the $600$-cell}
The 600-cell has the $f$-vector $(120,720,1200,600)$, by duality the
120-cell has the $f$-vector $(600,1200,720,120)$. 
Any 1-Hamiltonian surface in the 600-cell must
have 120 vertices, 720 edges and, consequently, 480 triangles (namely,
two out of five),
so it has Euler characteristic $\chi = -120$ and genus $g=61$.
We obtain the same genus in the 120-cell by counting 600 vertices, 1200 edges
and 480 pentagons (namely, two out of three). 
The same Euler characteristic would hold for a pinched surface if
there is any. We remark that similarly the 4-cube
admits a Hamiltonian surface of the same genus (namely, $g=1$) as
the 4-dimensional cross polytope. 

\subparagraph{Proposition 2}
{\it There is no $1$-Hamiltonian surface in the $2$-skeleton of the
  $120$-cell.
There is no pinched surface either since the vertex link of the $120$-cell
is too small for containing two disjoint circuits.}

\bigskip
The proof is a fairly simple procedure: In each vertex link of type $\{3,3\}$
the Hamiltonian
surface appears as a Hamiltonian circuit of length 4. This
is unique, up to symmetries of the tetrahedron and of the 120-cell
itself.
Note that two consecutive edges determine the circuit completely.
So without loss of generality we can start with such a unique vertex
link of the surface. This means we start with four pentagons
covering the star of one vertex.
In each of the four neighboring vertices this determines two
consecutive edges of the link there. It follows
that these circuits are uniquely determined as well
and that we can extend the beginning part of our surface,
now covering the stars of five vertices.
Successively this leads to a construction of such a surface.
However, after a few steps it ends at a contradiction.
Consequently, such a Hamiltonian surface does not exist.

\subparagraph{Proposition 3}
{\it There is no $1$-Hamiltonian surface in the $2$-skeleton of the
  $600$-cell.}

\bigskip
This proof is more involved since it 
uses the classification of all 17 distinct Hamiltonian
circuits in the icosahedron, up to symmetries
of it \cite[pp.\ 277 ff.]{H}.
If there is such a 1-Hamiltonian surface, then the link of each vertex
in it must be a Hamiltonian cycle in the vertex link of the 600-cell
which is an icosahedron.
We just have to see how these can fit together.
Starting with one arbitrary link one can try to extend the
triangulation to the neighbors. For the neighbors there are
forbidden 2-faces which has a consequence for the possible
types among the 17 for them.
After an exhaustive computer search it turned out 
that there is no way to fit all vertex links together.
Therefore such a surface does not exist.
At this point it must be left open whether there are 
1-Hamiltonian pinched surfaces in the 600-cell.
The reason is that there are too many possibilities for a splitting
into two, three or four cycles in the vertex link.
For a systematic search one would have to classify 
all these possibilities first.

\medskip
The GAP programs used for the algorithmic proof of Propositions 1, 2, 3
and details of the calculations
are available from the first author upon request.
%%%%%%%%%%%%%%%%%%%%%%SECTION 3  %%%%%%%%%%%%%%%%%%%%%%%%%%%%%%%%%%%%%%%%%%%
\subsection*{3. Hamiltonian submanifolds of cross polytopes}
The $d$-dimensional cross polytope $\beta^d$ (or the $d$-octahedron) is defined
as the convex hull of the $2d$ points
$$(0,\ldots,0,\pm1,0,\ldots,0) \in \mathbb{R}^d.$$
It is a simplicial and regular polytope, and it is
centrally-symmetric with $d$
diagonals, each between two antipodal points of type 
$(0,\ldots,0,1,0,\ldots,0)$ and $(0,\ldots,0,-1,0,\ldots,0).$
Its edge graph is the complete $d$-partite graph with two vertices in 
each partition, sometimes denoted by $K_2 * \cdots * K_2$.
See \cite {MMS} for properties of regular polytopes in general.
The $f$-vector of the cross polytope satisfies the equality
$$f_i(\beta^d) = 2^{i+1}{d \choose {i+1}}.$$
Consequently, any 1-Hamiltonian 2-manifold must have the following
beginning part of the $f$-vector:
$$f_0 = 2d, \ f_1 = 2d(d-1)$$
It follows that the Euler characteristic $\chi$ of the 2-manifold satisfies
$$2-\chi = 2-2d+2d(d-1)- \frac{4}{3}d(d-1) = \frac{2}{3}(d-1)(d-3).$$ 
These are the regular cases investigated in \cite{JR}.
In terms of the genus $g = \frac{1}{2}(2 - \chi)$ of an orientable surface 
this equation reads as
$$g = \frac{d-1}{1} \cdot \frac{d-3}{3}.$$ 
This remains valid for non-orientable surfaces 
if we assign the genus $\frac{1}{2}$ to the real projective
plane.
In any case $\chi$ can be an integer only if $d \equiv 0,1 (3)$.
The first possibilities, where all cases are actually realized 
by triangulations of closed orientable surfaces \cite{JR}, 
are indicated in Table 3.

\begin{center}{\small
\begin{tabular}{|r||r|r|}\hline
$d$ & $2 - \chi$ & genus $g$\\ \hline

3 & 0 & 0\\
4 & 2 & 1 \\
6 & 10 & 5\\
7 & 16 & 8\\
9 & 32 & 16\\
10 & 42 & $3 \cdot 7 = 21$\\
12 & 66 & $3 \cdot 11 = 33$\\
13 & 80 & $8 \cdot 5  = 40$\\
15 & 112 & $8 \cdot 7  = 56$\\
16 & 120 & $4 \cdot 3 \cdot 5 = 60$\\
18 & 170 & $5 \cdot 17 = 85$\\
19 & 192 & $32 \cdot 3  = 96$\\
21 & 240 & $8 \cdot 3 \cdot 5 = 120$\\
22 & 266 & $7 \cdot 19 = 133$\\

\hline
\end{tabular}}

\medskip
Table 3: Regular cases of 1-Hamiltonian 2-manifolds
\end{center}

 \medskip
Similarly, any 2-Hamiltonian 4-manifold must have the following
beginning part of the $f$-vector:
$$f_0 = 2d, \ f_1 = 2d(d-1), f_2 = \frac{4}{3}d(d-1)(d-2)$$
It follows that the Euler characteristic $\chi$ satisfies
$$10(\chi - 2) = f_2 - 4f_1 + 10f_0 - 20 = \frac{4}{3}d(d-1)(d-2)
-8d(d-1) + 20d - 20 =  \frac{4}{3}(d-1)(d-3)(d-5).$$
If we introduce the ``genus'' $g = \frac{1}{2}(\chi-2)$ 
of a simply connected 4-manifold as the number
of copies of $S^2 \times S^2$ which are necessary to form a connected
sum with Euler characteristic $\chi$, then this equation reads as
$$g = \frac{d-1}{1} \cdot \frac{d-3}{3} \cdot \frac{d-5}{5}.$$ 
These are the ``regular cases''.
Again the complex projective plane would have genus $\frac{1}{2}$ here.
Recall that any 2-Hamiltonian 4-manifold in the boundary of a convex
polytope is simply connected since the 2-skeleton is. Therefore the 
``genus'' equals half of the second Betti number.

\bigskip
Moreover, there is an Upper Bound Theorem and a Lower Bound Theorem
as follows:

\subparagraph{Theorem 1} (E.~Sparla \cite{Sp1})

{\it If a triangulation of a $4$-manifold occurs as a $2$-Hamiltonian
subcomplex of a centrally-symmetric
simplicial $d$-polytope then the following inequality holds
$$\frac{1}{2}\big(\chi(M) - 2\big) \geq \frac{d-1}{1} \cdot \frac{d-3}{3} \cdot \frac{d-5}{5}.$$ 
Moreover, for $d \geq 6$ equality is possible only if the polytope
is affinely equivalent to the $d$-dimensional cross polytope.

\medskip
If there is a triangulation of a $4$-manifold with a fixed point free
involution then the number $n$ of vertices is even, i.e., $n = 2d$, and
the opposite inequality holds
$$\frac{1}{2}\big(\chi(M) - 2\big) \leq \frac{d-1}{1} \cdot \frac{d-3}{3} \cdot \frac{d-5}{5}.$$ 
Moreover, equality in this inequality implies that the manifold
can be regarded as a $2$-Hamiltonian subcomplex of the $d$-dimensional
cross polytope.}

\bigskip
Remark. The case of equality in either of these inequalities corresponds
to the ``regular cases''. Sparla's original equation
$4^3{{\frac{1}{2}(d-1)} \choose {3}} = 10(\chi(M) - 2)$
is equivalent to the one above.

\bigskip
By analogy, any $k$-Hamiltonian $2k$-manifold in the $d$-dimensional
cross polytope satisfies the equation
$$(-1)^k\frac{1}{2}\big(\chi - 2\big) = \frac{d-1}{1} \cdot \frac{d-3}{3} \cdot \frac{d-5}{5} \cdot
\ \cdots \ \cdot \frac{d-2k-1}{2k+1}.$$ 
It is necessarily $(k-1)$-connected which implies that the left
hand side is half of the middle Betti number which 
is nothing but the ``genus''.
Furthermore, there is a conjectured Upper Bound Theorem and a Lower Bound
Theorem generalizing Theorem 1 where the inequality has to be replaced by
$$(-1)^k\frac{1}{2}\big(\chi - 2\big) \geq\frac{d-1}{1} \cdot \frac{d-3}{3} \cdot \frac{d-5}{5} \cdot
\ \cdots \ \cdot \frac{d-2k-1}{2k+1}$$ or 
$$(-1)^k\frac{1}{2}\big(\chi - 2\big) \leq \frac{d-1}{1} \cdot \frac{d-3}{3} \cdot \frac{d-5}{5} \cdot
\ \cdots \ \cdot \frac{d-2k-1}{2k+1},$$ 
respectively, see \cite{Sp2}, \cite{N2}.
The discussion of the cases of equality is exactly the same.
Sparla's original version
$$4^{k+1}{{\frac{1}{2}(d-1)} \choose {k+1}} = {{2k+1} \choose {k+1}}(-1)^k\big(\chi(M) - 2\big)$$
is equivalent to the one above.
In particular, for any $k$ one of the ``regular cases'' 
is the case of a sphere product
$S^k \times S^k$ with $(-1)^k(\chi - 2) = 2$ 
(or ``genus'' $g=1$) and $d = 2k+2$.
So far examples are available for $1 \leq k \leq 4$, even with a vertex
transitive automorphism group
see \cite{Lu}, \cite{K-Lu}.
We hope that for $k \geq 5$ there will be similar examples as well,
compare Section 6.
%%%%%%%%%%%%%%%%%%%%%% SECTION 4 %%%%%%%%%%%%%%%%%%%%%%%%%%%%%%%%%%%%%%%%
\subsection*{4. $2$-Hamiltonian $4$-manifolds in cross polytopes}
In the case of 2-Hamiltonian 4-manifolds as subcomplexes of
the $d$-dimensional cross polytope we have the ``regular cases''
of equality 
$g =\frac{1}{2}(\chi -2)= \frac{d-1}{1} \cdot \frac{d-3}{3} \cdot \frac{d-5}{5}.$ 
Here $\chi$ can be an integer only if $d \equiv 0,1,3 (5)$.
Table 4 indicates the first possibilities:

\begin{center}{\small
\begin{tabular}{|r||r|r|c|}\hline
$d$ & $\chi - 2$ & ``genus'' $g$&existence\\ \hline

5 & 0 & 0&$S^4 = \partial \beta^5$\\
6 & 2 & 1&$S^2 \times S^2$ \cite{Sp1},\cite{La-Sp}\\
8 & 14 & 7&new (Thm.\ 2)\\
10 & 42 &$3 \cdot 7 = 21$&see Remark 2\\
11 & 64 & 32&?\\
13 & 128 & 64 &?\\
15 & 224 & $16 \cdot 7 = 112$&?\\
16 & 286 & $11 \cdot 13 = 143$&?\\
18 & 442 & $13 \cdot 17 = 221$&?\\
20 & 646 & $17 \cdot 19 = 323$&?\\
21 & 720 & $8 \cdot 5 \cdot 9= 360$&?\\
23 & 1056 & $16 \cdot 3 \cdot 11 = 528$&?\\
25 & 1408 & $64 \cdot 11 = 704$&?\\
26 & 1610 & $5 \cdot 7 \cdot 23 = 805$&?\\
28 & 2070 & $5 \cdot 9 \cdot 23 = 1035$&?\\
30 & 2610 & $5 \cdot 9 \cdot 29 = 1305$&?\\
%31 & 2800 & $8 \cdot 7 \cdot 25= 1400$&?\\
%33 & 3584 & $256 \cdot 7= 1792$&?\\
%35 & 4352 & $128 \cdot 17 = 2176$&?\\
%36 & 4774 & $7 \cdot 11 \cdot 31 = 2387$&?\\
%38 & 5698 & $7 \cdot 11 \cdot 37 = 2849$&?\\
%40 & 6734 & $7 \cdot 13 \cdot 37 = 3367$&?\\
\hline
\end{tabular}}

\medskip
Table 4: Regular cases of 2-Hamiltonian 4-manifolds
\end{center}

\subparagraph{Theorem 2}
{\it There is a $16$-vertex triangulation of a $4$-manifold $M \cong
  (S^2 \times S^2)^{\#7}$ which can be
regarded as a centrally-symmetric and $2$-Hamiltonian 
subcomplex of the $8$-dimensional cross polytope. 
As one of the ``regular cases'' it 
satisfies equality in Sparla's inequalities in Theorem $1$ with 
the ``genus'' $g = 7$ and with $d=8$. }

\bigskip 
{\sc Proof.} 
Any 2-Hamiltonian subcomplex of a convex polytope is simply connected
\cite[3.8]{Ku1}. Therefore such an $M$, if it exists, must
be simply connected, in particular $H_1(M) = H_3(M) = 0$.
In accordance with Sparla's inequalities, the Euler characteristic 
$\chi(M) = 16$ tells us that
the middle homology group is $H_2(M,\mathbb{Z}) \cong \mathbb{Z}^{14}$.
The topological type of $M$ is then uniquely determined by the
intersection form.
If the intersection form is even then by Rohlin's theorem 
the signature must be zero, which implies 
that $M$ is homeomorphic to 
the connected sum of $7$ copies of $S^2 \times S^2$,
see \cite{Sa}.
If the intersection form is odd then $M$ is a connected sum
of 14 copies of $\pm\mathbb{C}P^2$.
We will show that the intersection form of our example is even.

\medskip
The induced polyhedral embedding of this manifold into 8-space is tight
since the intersection with any open halfspace is connected and simply
connected, compare Section 5 below.
No smooth tight embedding of this manifold 
into 8-space can exist, see \cite{Th}.
Consequently, this embedding of $M$ into 8-space is smoothable as far
as the PL structure is concerned but it is not tightly smoothable.

\bigskip  
The $f$-vector $f = (16, 112, 448, 560, 224)$ of this example is
uniquely determined already by the requirement of
16 vertices and the condition to be 2-Hamiltonian
in the 8-dimensional cross polytope.
In particular there are 8 missing edges corresponding to the 8 diagonals 
of the cross polytope which are pairwise disjoint.

\medskip
Assuming a vertex-transitive automorphism group,
the example was found by using the software of F.~H.~Lutz described
in \cite{Lu}. 
The combinatorial automorphism group $G$ of our example is of order 128.
With this particular automorphism group the example is unique.
The special element 
$$\zeta = (1  \ 2)(3  \ 4)(5  \ 6)(7  \ 8)(9  \ 10)(11  \ 12)(13  \ 14)(15  \ 16)$$
acts on $M$ without fixed points. It interchanges the endpoints
of each diagonal and, therefore, can be regarded
as the antipodal mapping sending each vertex of the 8-dimensional
cross polytope to its antipodal vertex in such a way that it is
compatible with the subcomplex $M$. 
A normal subgroup $H$ isomorphic to 
$C_2 \oplus C_2 \oplus C_2 \oplus C_2$ 
acts simply transitively on the 16 vertices. 
The isotropy group $G_0$ fixing
one vertex (and, simultaneously, its antipodal vertex) 
is isomorphic to the dihedral group of order 8.
The group itself is a semidirect product between $H$ and $G_0$.
In more detail the example is given by the three $G$-orbits of the 4-simplices 
$$\langle 1\, 3\, 5\, 7\, 9\rangle_{128},\ \langle 1\, 3\, 5\, 9\, 13\rangle_{64},\ \langle 1\, 3\, 5\, 7\,15\rangle_{32}$$
with altogether $128 + 64 + 32 = 224$ simplices, 
each given by a 5-tuple of vertices out of $\{1,2,3, \ldots, 15,16\}$.
The group $G \cong ((((C_4 \oplus C_2):C_2):C_2):C_2):C_2$
of order 128 is generated by the three permutations
$\alpha = (1 \ 12 \ 16 \ 14  \ 2 \ 11 \ 15 \ 13)(3 \ 10 \ 6  \ 8 \ 4 \ 9 \ 5  \ 7),$ 
$\beta = (1 \ 6 \ 2 \ 5)(7   \ 9 \ 2 \ 14)(8 \ 10 \ 11 \ 13)(15 \ 16),$
$\gamma = (1  \ 12 \ 3 \ 14)(2 \ 11 \ 4 \ 13)(5 \ 7 \ 16 \ 10)(6  \ 8 \ 15 \ 9).$

\bigskip
The complete list of all 224 top-dimensional simplices is the following:

\bigskip
\begin{center}
{\scriptsize
\begin{tabular}{@{}l@{}l@{}l@{}l@{}l@{}l@{}}
{$\langle 1\,3\,5\,7\,9 \rangle$},&{$\langle 1\,3\,5\,7\,15 \rangle$},&{$\langle 1\,3\,5\,8\,13 \rangle$},&{$\langle 1\,3\,5\,8\,15 \rangle$},&{$\langle 1\,3\,5\,9\,13 \rangle$},&{$\langle 1\,3\,6\,8\,10 \rangle$},\\
{$\langle 1\,3\,6\,8\,12 \rangle$},&{$\langle 1\,3\,6\,9\,12 \rangle$},&{$\langle 1\,3\,6\,9\,16 \rangle$},&{$\langle 1\,3\,6\,10\,16 \rangle$},&{$\langle 1\,3\,7\,9\,15 \rangle$},&{$\langle 1\,3\,8\,10\,16 \rangle$},\\
{$\langle 1\,3\,8\,11\,14 \rangle$},&{$\langle 1\,3\,8\,11\,16 \rangle$},&{$\langle 1\,3\,8\,12\,13 \rangle$},&{$\langle 1\,3\,8\,14\,15 \rangle$},&{$\langle 1\,3\,9\,11\,14 \rangle$},&{$\langle 1\,3\,9\,11\,16 \rangle$},\\
{$\langle 1\,3\,9\,12\,13 \rangle$},&{$\langle 1\,3\,9\,14\,15 \rangle$},&  {$\langle 1\,4\,5\,9\,12 \rangle$},&{$\langle 1\,4\,5\,9\,13 \rangle$},&{$\langle 1\,4\,5\,11\,13 \rangle$},&{$\langle 1\,4\,5\,11\,16 \rangle$},\\
{$\langle 1\,4\,5\,12\,16 \rangle$},&  {$\langle 1\,4\,6\,8\,12 \rangle$},&{$\langle 1\,4\,6\,8\,13 \rangle$},&{$\langle 1\,4\,6\,12\,14 \rangle$},&{$\langle 1\,4\,6\,13\,15 \rangle$},&{$\langle 1\,4\,6\,14\,15 \rangle$},\\
{$\langle 1\,4\,7\,10\,12 \rangle$},&{$\langle 1\,4\,7\,10\,13 \rangle$},&{$\langle 1\,4\,7\,12\,15 \rangle$},&{$\langle 1\,4\,7\,13\,15 \rangle$},&{$\langle 1\,4\,8\,9\,12 \rangle$},&{$\langle 1\,4\,8\,9\,13 \rangle$},\\
{$\langle 1\,4\,10\,12\,16 \rangle$},&{$\langle 1\,4\,10\,13\,16 \rangle$},&{$\langle 1\,4\,11\,13\,16 \rangle$},&{$\langle 1\,4\,12\,14\,15 \rangle$},&{$\langle 1\,5\,7\,9\,12 \rangle$},&{$\langle 1\,5\,7\,10\,12 \rangle$},\\
{$\langle 1\,5\,7\,10\,15 \rangle$},&{$\langle 1\,5\,8\,11\,13 \rangle$},&{$\langle 1\,5\,8\,11\,14 \rangle$},&{$\langle 1\,5\,8\,14\,15 \rangle$},&{$\langle 1\,5\,10\,12\,16 \rangle$},&{$\langle 1\,5\,10\,14\,15 \rangle$},\\
{$\langle 1\,5\,10\,14\,16 \rangle$},&{$\langle 1\,5\,11\,14\,16 \rangle$},&{$\langle 1\,6\,7\,10\,13 \rangle$},&{$\langle 1\,6\,7\,10\,16 \rangle$},&{$\langle 1\,6\,7\,11\,15 \rangle$},&{$\langle 1\,6\,7\,11\,16 \rangle$},\\
{$\langle 1\,6\,7\,13\,15 \rangle$},&{$\langle 1\,6\,8\,10\,13 \rangle$},&{$\langle 1\,6\,9\,11\,14 \rangle$},&{$\langle 1\,6\,9\,11\,16 \rangle$},&{$\langle 1\,6\,9\,12\,14 \rangle$},&{$\langle 1\,6\,11\,14\,15 \rangle$},\\
{$\langle 1\,7\,9\,12\,15 \rangle$},&{$\langle 1\,7\,10\,11\,14 \rangle$},&{$\langle 1\,7\,10\,11\,15 \rangle$},&{$\langle 1\,7\,10\,14\,16 \rangle$},&{$\langle 1\,7\,11\,14\,16 \rangle$},&{$\langle 1\,8\,9\,12\,13 \rangle$},\\
{$\langle 1\,8\,10\,13\,16 \rangle$},&{$\langle 1\,8\,11\,13\,16 \rangle$},&{$\langle 1\,9\,12\,14\,15 \rangle$},&{$\langle 1\,10\,11\,14\,15 \rangle$},&{$\langle 2\,3\,5\,7\,11 \rangle$},&{$\langle 2\,3\,5\,7\,14 \rangle$},\\
{$\langle 2\,3\,5\,11\,13 \rangle$},&{$\langle 2\,3\,5\,13\,16 \rangle$},&{$\langle 2\,3\,5\,14\,16 \rangle$},&{$\langle 2\,3\,6\,10\,11 \rangle$},&{$\langle 2\,3\,6\,10\,14 \rangle$},&{$\langle 2\,3\,6\,11\,15 \rangle$},\\
{$\langle 2\,3\,6\,12\,14 \rangle$},&{$\langle 2\,3\,6\,12\,15 \rangle$},&{$\langle 2\,3\,7\,10\,11 \rangle$},&{$\langle 2\,3\,7\,10\,14 \rangle$},&{$\langle 2\,3\,8\,9\,11 \rangle$},&{$\langle 2\,3\,8\,9\,14 \rangle$},\\
{$\langle 2\,3\,8\,11\,16 \rangle$},&{$\langle 2\,3\,8\,14\,16 \rangle$},&{$\langle 2\,3\,9\,11\,15 \rangle$},&{$\langle 2\,3\,9\,14\,15 \rangle$},&{$\langle 2\,3\,11\,13\,16 \rangle$},&{$\langle 2\,3\,12\,14\,15 \rangle$},\\
{$\langle 2\,4\,5\,7\,9 \rangle$},&{$\langle 2\,4\,5\,7\,11 \rangle$},&{$\langle 2\,4\,5\,9\,15 \rangle$},&{$\langle 2\,4\,5\,10\,11 \rangle$},&{$\langle 2\,4\,5\,10\,15 \rangle$},&{$\langle 2\,4\,6\,7\,14 \rangle$},\\
{$\langle 2\,4\,6\,7\,16 \rangle$},&{$\langle 2\,4\,6\,8\,10 \rangle$},&{$\langle 2\,4\,6\,8\,16 \rangle$},&{$\langle 2\,4\,6\,10\,14 \rangle$},&{$\langle 2\,4\,7\,9\,15 \rangle$},&{$\langle 2\,4\,7\,11\,14 \rangle$},\\
{$\langle 2\,4\,7\,12\,13 \rangle$},&{$\langle 2\,4\,7\,12\,15 \rangle$},&{$\langle 2\,4\,7\,13\,16 \rangle$},&{$\langle 2\,4\,8\,10\,16 \rangle$},&{$\langle 2\,4\,10\,11\,14 \rangle$},&{$\langle 2\,4\,10\,12\,13 \rangle$},\\
{$\langle 2\,4\,10\,12\,15 \rangle$},&{$\langle 2\,4\,10\,13\,16 \rangle$},&{$\langle 2\,5\,7\,9\,14 \rangle$},&{$\langle 2\,5\,8\,9\,14 \rangle$},&{$\langle 2\,5\,8\,9\,15 \rangle$},&{$\langle 2\,5\,8\,12\,15 \rangle$},\\
{$\langle 2\,5\,8\,12\,16 \rangle$},&{$\langle 2\,5\,8\,14\,16 \rangle$},&{$\langle 2\,5\,10\,11\,13 \rangle$},&{$\langle 2\,5\,10\,12\,13 \rangle$},&{$\langle 2\,5\,10\,12\,15 \rangle$},&{$\langle 2\,5\,12\,13\,16 \rangle$},\\
{$\langle 2\,6\,7\,12\,13 \rangle$},&{$\langle 2\,6\,7\,12\,14 \rangle$},&{$\langle 2\,6\,7\,13\,16 \rangle$},&{$\langle 2\,6\,8\,9\,11 \rangle$},&{$\langle 2\,6\,8\,9\,16 \rangle$},&{$\langle 2\,6\,8\,10\,11 \rangle$},\\
{$\langle 2\,6\,9\,11\,15 \rangle$},&{$\langle 2\,6\,9\,13\,15 \rangle$},&{$\langle 2\,6\,9\,13\,16 \rangle$},&{$\langle 2\,6\,12\,13\,15 \rangle$},&{$\langle 2\,7\,9\,14\,15 \rangle$},&{$\langle 2\,7\,10\,11\,14 \rangle$},\\
{$\langle 2\,7\,12\,14\,15 \rangle$},&{$\langle 2\,8\,9\,12\,13 \rangle$},&{$\langle 2\,8\,9\,12\,16 \rangle$},&{$\langle 2\,8\,9\,13\,15 \rangle$},&{$\langle 2\,8\,10\,11\,16 \rangle$},&{$\langle 2\,8\,12\,13\,15 \rangle$},\\
{$\langle 2\,9\,12\,13\,16 \rangle$},&{$\langle 2\,10\,11\,13\,16 \rangle$},&{$\langle 3\,5\,7\,9\,11 \rangle$},&{$\langle 3\,5\,7\,10\,12 \rangle$},&{$\langle 3\,5\,7\,10\,15 \rangle$},&{$\langle 3\,5\,7\,12\,16 \rangle$},\\
{$\langle 3\,5\,7\,14\,16 \rangle$},&{$\langle 3\,5\,8\,13\,15 \rangle$},&{$\langle 3\,5\,9\,11\,13 \rangle$},&{$\langle 3\,5\,10\,12\,13 \rangle$},&{$\langle 3\,5\,10\,13\,15 \rangle$},&{$\langle 3\,5\,12\,13\,16 \rangle$},\\
{$\langle 3\,6\,7\,10\,13 \rangle$},&{$\langle 3\,6\,7\,10\,16 \rangle$},&{$\langle 3\,6\,7\,12\,13 \rangle$},&{$\langle 3\,6\,7\,12\,16 \rangle$},&{$\langle 3\,6\,8\,10\,14 \rangle$},&{$\langle 3\,6\,8\,12\,14 \rangle$},\\
{$\langle 3\,6\,9\,12\,16 \rangle$},&{$\langle 3\,6\,10\,11\,15 \rangle$},&{$\langle 3\,6\,10\,13\,15 \rangle$},&{$\langle 3\,6\,12\,13\,15 \rangle$},&{$\langle 3\,7\,9\,11\,15 \rangle$},&{$\langle 3\,7\,10\,11\,15 \rangle$},\\
{$\langle 3\,7\,10\,12\,13 \rangle$},&{$\langle 3\,7\,10\,14\,16 \rangle$},&{$\langle 3\,8\,9\,11\,14 \rangle$},&{$\langle 3\,8\,10\,14\,16 \rangle$},&{$\langle 3\,8\,12\,13\,15 \rangle$},&{$\langle 3\,8\,12\,14\,15 \rangle$},\\
{$\langle 3\,9\,11\,13\,16 \rangle$},&{$\langle 3\,9\,12\,13\,16 \rangle$},&{$\langle 4\,5\,7\,9\,13 \rangle$},&{$\langle 4\,5\,7\,11\,13 \rangle$},&{$\langle 4\,5\,8\,9\,14 \rangle$},&{$\langle 4\,5\,8\,9\,15 \rangle$},\\
{$\langle 4\,5\,8\,11\,14 \rangle$},&{$\langle 4\,5\,8\,11\,15 \rangle$},&{$\langle 4\,5\,9\,12\,16 \rangle$},&{$\langle 4\,5\,9\,14\,16 \rangle$},&{$\langle 4\,5\,10\,11\,15 \rangle$},&{$\langle 4\,5\,11\,14\,16 \rangle$},\\
{$\langle 4\,6\,7\,14\,16 \rangle$},&{$\langle 4\,6\,8\,9\,11 \rangle$},&{$\langle 4\,6\,8\,9\,16 \rangle$},&{$\langle 4\,6\,8\,10\,12 \rangle$},&{$\langle 4\,6\,8\,11\,15 \rangle$},&{$\langle 4\,6\,8\,13\,15 \rangle$},\\
{$\langle 4\,6\,9\,11\,14 \rangle$},&{$\langle 4\,6\,9\,14\,16 \rangle$},&{$\langle 4\,6\,10\,12\,14 \rangle$},&{$\langle 4\,6\,11\,14\,15 \rangle$},&{$\langle 4\,7\,9\,13\,15 \rangle$},&{$\langle 4\,7\,10\,12\,13 \rangle$},\\
{$\langle 4\,7\,11\,13\,16 \rangle$},&{$\langle 4\,7\,11\,14\,16 \rangle$},&{$\langle 4\,8\,9\,11\,14 \rangle$},&{$\langle 4\,8\,9\,12\,16 \rangle$},&{$\langle 4\,8\,9\,13\,15 \rangle$},&{$\langle 4\,8\,10\,12\,16 \rangle$},\\
{$\langle 4\,10\,11\,14\,15 \rangle$},&{$\langle 4\,10\,12\,14\,15 \rangle$},&{$\langle 5\,7\,9\,11\,13 \rangle$},&{$\langle 5\,7\,9\,12\,16 \rangle$},&{$\langle 5\,7\,9\,14\,16 \rangle$},&{$\langle 5\,8\,10\,12\,14 \rangle$},\\
{$\langle 5\,8\,10\,12\,16 \rangle$},&{$\langle 5\,8\,10\,14\,16 \rangle$},&{$\langle 5\,8\,11\,13\,15 \rangle$},&{$\langle 5\,8\,12\,14\,15 \rangle$},&{$\langle 5\,10\,11\,13\,15 \rangle$},&{$\langle 5\,10\,12\,14\,15 \rangle$},\\
{$\langle 6\,7\,9\,11\,13 \rangle$},&{$\langle 6\,7\,9\,11\,15 \rangle$},&{$\langle 6\,7\,9\,13\,15 \rangle$},&{$\langle 6\,7\,11\,13\,16 \rangle$},&{$\langle 6\,7\,12\,14\,16 \rangle$},&{$\langle 6\,8\,10\,11\,15 \rangle$},\\
{$\langle 6\,8\,10\,12\,14 \rangle$},&{$\langle 6\,8\,10\,13\,15 \rangle$},&{$\langle 6\,9\,11\,13\,16 \rangle$},&{$\langle 6\,9\,12\,14\,16 \rangle$},&{$\langle 7\,9\,12\,14\,15 \rangle$},&{$\langle 7\,9\,12\,14\,16 \rangle$},\\
{$\langle 8\,10\,11\,13\,15 \rangle$},&{$\langle 8\,10\,11\,13\,16 \rangle$}.&&&&
\end{tabular}
}
\end{center}

The link of the vertex 16 is the following simplicial 3 sphere
with 70 tetrahedra:
\begin{center}
{\scriptsize
\begin{tabular}{@{}l@{}l@{}l@{}l@{}l@{}l@{}l@{}l@{}}
{$\langle 1\,3\,6\,9 \rangle$},&{$\langle 1\,3\,6\,10 \rangle$},&{$\langle 1\,3\,8\,10 \rangle$},&{$\langle 1\,3\,8\,11 \rangle$},&{$\langle 1\,3\,9\,11 \rangle$},&{$\langle 1\,4\,5\,11 \rangle$},&{$\langle 1\,4\,5\,12 \rangle$},&{$\langle 1\,4\,10\,12 \rangle$},\\
{$\langle 1\,4\,10\,13 \rangle$},&{$\langle 1\,4\,11\,13 \rangle$},&{$\langle 1\,5\,10\,12 \rangle$},&{$\langle 1\,5\,10\,14 \rangle$},&{$\langle 1\,5\,11\,14 \rangle$},&{$\langle 1\,6\,7\,10 \rangle$},&{$\langle 1\,6\,7\,11 \rangle$},&{$\langle 1\,6\,9\,11 \rangle$},\\
{$\langle 1\,7\,10\,14 \rangle$},&{$\langle 1\,7\,11\,14 \rangle$},&{$\langle 1\,8\,10\,13 \rangle$},&{$\langle 1\,8\,11\,13 \rangle$},&{$\langle 2\,3\,5\,13 \rangle$},&{$\langle 2\,3\,5\,14 \rangle$},&{$\langle 2\,3\,8\,11 \rangle$},&{$\langle 2\,3\,8\,14 \rangle$},\\
{$\langle 2\,3\,11\,13 \rangle$},&{$\langle 2\,4\,6\,7 \rangle$},&{$\langle 2\,4\,6\,8 \rangle$},&{$\langle 2\,4\,7\,13 \rangle$},&{$\langle 2\,4\,8\,10 \rangle$},&{$\langle 2\,4\,10\,13 \rangle$},&{$\langle 2\,5\,8\,12 \rangle$},&{$\langle 2\,5\,8\,14 \rangle$},\\
{$\langle 2\,5\,12\,13 \rangle$},&{$\langle 2\,6\,7\,13 \rangle$},&{$\langle 2\,6\,8\,9 \rangle$},&{$\langle 2\,6\,9\,13 \rangle$},&{$\langle 2\,8\,9\,12 \rangle$},&{$\langle 2\,8\,10\,11 \rangle$},&{$\langle 2\,9\,12\,13 \rangle$},&{$\langle 2\,10\,11\,13 \rangle$},\\
{$\langle 3\,5\,7\,12 \rangle$},&{$\langle 3\,5\,7\,14 \rangle$},&{$\langle 3\,5\,12\,13 \rangle$},&{$\langle 3\,6\,7\,10 \rangle$},&{$\langle 3\,6\,7\,12 \rangle$},&{$\langle 3\,6\,9\,12 \rangle$},&{$\langle 3\,7\,10\,14 \rangle$},&{$\langle 3\,8\,10\,14 \rangle$},\\
{$\langle 3\,9\,11\,13 \rangle$},&{$\langle 3\,9\,12\,13 \rangle$},&{$\langle 4\,5\,9\,12 \rangle$},&{$\langle 4\,5\,9\,14 \rangle$},&{$\langle 4\,5\,11\,14 \rangle$},&{$\langle 4\,6\,7\,14 \rangle$},&{$\langle 4\,6\,8\,9 \rangle$},&{$\langle 4\,6\,9\,14 \rangle$},\\
{$\langle 4\,7\,11\,13 \rangle$},&{$\langle 4\,7\,11\,14 \rangle$},&{$\langle 4\,8\,9\,12 \rangle$},&{$\langle 4\,8\,10\,12 \rangle$},&{$\langle 5\,7\,9\,12 \rangle$},&{$\langle 5\,7\,9\,14 \rangle$},&{$\langle 5\,8\,10\,12 \rangle$},&{$\langle 5\,8\,10\,14\rangle$},\\
{$\langle 6\,7\,11\,13 \rangle$},&{$\langle 6\,7\,12\,14 \rangle$},&{$\langle 6\,9\,11\,13 \rangle$},&{$\langle 6\,9\,12\,14 \rangle$},&{$\langle 7\,9\,12\,14 \rangle$},&{$\langle 8\,10\,11\,13\rangle$}.&&
\end{tabular}
}
\end{center}

It remains to prove two facts:

\medskip
\underline{Claim 1.} The link of the vertex 16 is a combinatorial 3-sphere.
This implies that $M$ is a PL-manifold
since all vertices are equivalent under the action of the automorphism group.

\medskip
A computer algorithm gave a positive answer: the link of the vertex 16
is combinatorially equivalent to the boundary of a 4-simplex by 
bistellar moves. This method is described in \cite{Bj-Lu} and \cite[1.3]{Lu}.

\medskip
\underline{Claim 2.} The intersection form of $M$ is even or, 
equivalently, the second
Stiefel-Whitney class of $M$ vanishes.
This implies that $M$ is homeomorphic to the connected sum of $7$ copies
of $S^2 \times S^2$.

\medskip
There is an algorithm for calculating the second 
Stiefel-Whitney class \cite{GT}.
There is also an computer algorithm implemented in {\tt polymake} \cite{G-J}, compare \cite{J} for determining the intersection form itself. 
The latter algorithm gave the following answer:
The intersection form of $M$ is even, and the signature is zero.
\hfill $\Box$

\medskip
In order to illustrate the intersection form on the second homology
we consider the link of the vertex 16, as given above. 
By the tightness condition special homology classes are
represented by the empty tetrahedra $c_1=\langle 7\,10\,11\,16 \rangle$
and $d_1 = \langle 8\,12\,13\,16 \rangle$ 
which are interchanged by the element
$$\delta = (1 \ 2)(5 \ 6)(7 \ 12)(8 \ 11)(9 \ 14)(10 \ 13)$$ 
of the automorphism group.
The intersection number of these two equals
the linking number of the empty triangles  $\langle 7\,10\,11 \rangle$
and $\langle 8\,12\,13 \rangle$ in the link of 16.
The two subsets in the link spanned by 
$1,5,7,10,11,14$ and $2,6,8,9,12,13$, respectively,
are homotopy circles interchanged by $\delta$. 
The intermediate subset
of points in the link of 16 which is invariant under $\delta$ 
is the torus depicted in Figure 4. The set of points which are fixed 
by $\delta$ are represented 
as the horizontal $(1,1)$-curve in this torus, the element $\delta$ itself
appears as the reflection along that fixed curve.
This torus shrinks down to the homotopy circle on either of the sides
which are spanned by $1,5,7,10,11,14$ and $2,6,8,9,12,13$, respectively. 
The empty triangles $\langle 7\,10\,11 \rangle$
and $\langle 8\,12\,13 \rangle$ also represent the same homotopy circles.
Since the link is a 3-sphere these two are linked with linking number $\pm 1$.
As a result we get for the intersection form
$c_1 \cdot d_1 = \pm 1$. These two empty tetrahedra $c_1$ and $d_1$ 
are not homologous
to one another in $M$. Each one can be perturbed into a disjoint position
such that the self linking number is zero: 
$c_1 \cdot c_1 = d_1 \cdot d_1 = 0$. Therefore $c_1,d_1$ represent
a part of the intersection form isomorphic with 
$\pm {{0 \ \ 1}\choose{1 \ \ 0}}$.
This situation is transferred to the intersection form of other
generators by the automorphism group. 
As a result we have seven copies of the matrix as a direct sum.

\begin{center}
\hspace*{-0.4cm}
\epsfig{file=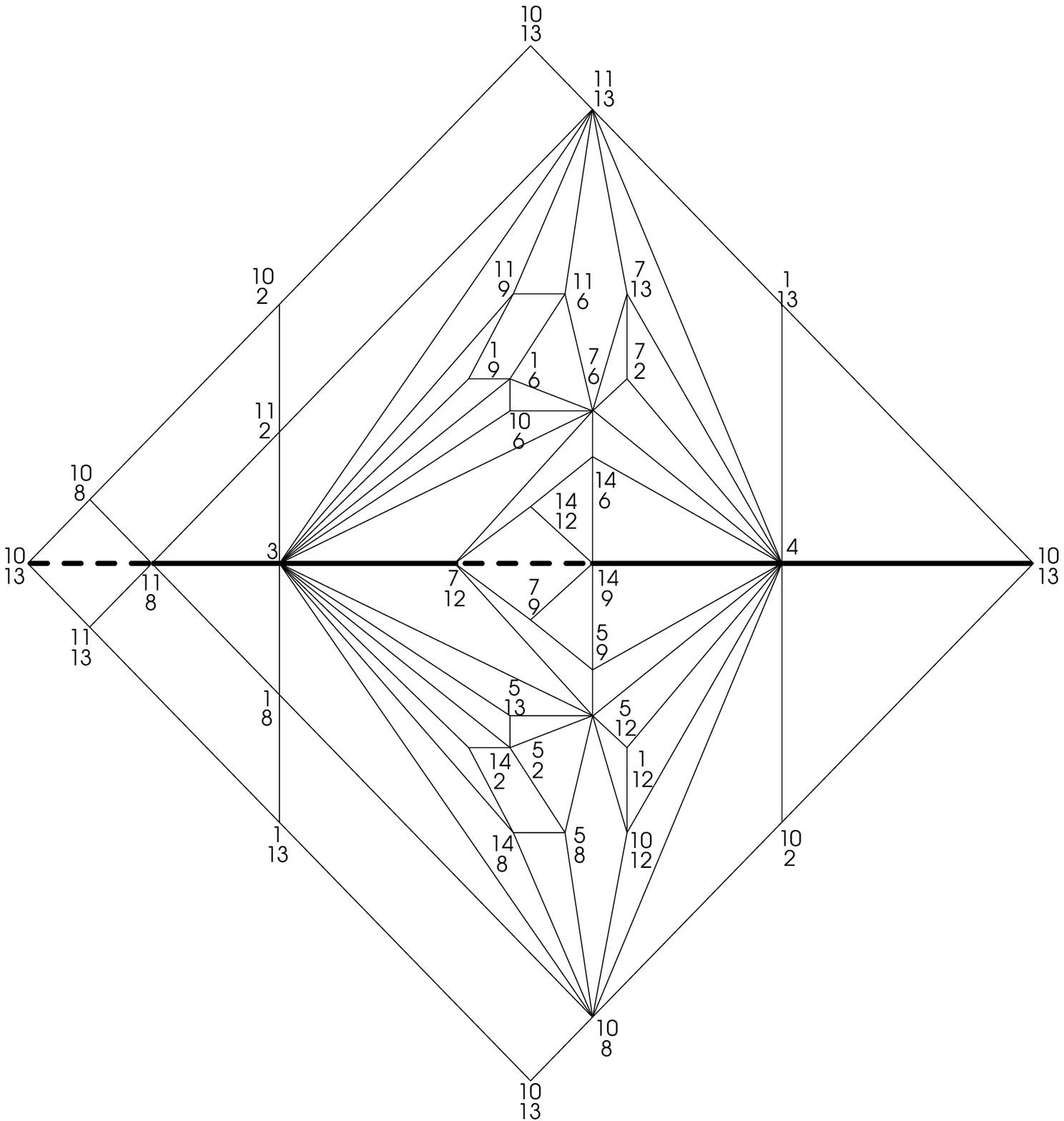,height=13.5cm}

\medskip
{Figure 4: The intermediate torus in the link of 16, invariant under
  the reflection $\delta$}
\end{center}

\bigskip
{\bf Remark 1:} 
Looking at the action of the automorphism group $G$
on the free abelian group $H_2(M,\mathbb{Z}) \cong \mathbb{Z}^{14}$ 
we get on the 17 conjugacy
classes of $G$ the following 
character values
  $$  (14,-2,-2,-2,2,-2,6,-2,-2,-2,6,0,0,0,0,0,0) \ .$$
Denote by $\chi $ the corresponding ordinary character.  
Using the character table\footnote{We thank Wolfgang Kimmerle for 
helpful comments concerning group representations.} of $G$ given by 
GAP \cite{GAP} and the
orthogonality relations this character decomposes into a sum
of five irreducible ordinary characters as follows
  $$ \chi = \chi_2 + \chi_3 + \chi_{13} + \chi_{14} + \chi_{17} \ .$$
This shows that $\mathbb{C} \otimes_{\mathbb{Z}} H_2(M,\mathbb{Z})$ 
is a cyclic 
$\mathbb{C} G$ - module.
It may be interesting to find a geometric explanation for this. 
 The involved irreducible characters are as
follows:
$$
\begin{array}{l|rrrrrrrrrrrrrrrrr}
  & 1a & 2a& 2b & 2c& 4a & 2d& 2e& 4b& 4c & 4d& 2f& 4e& 4f & 4g & 4h &
  8a & 2g \\ \hline
\chi_2 & 1 & -1 & 1 & -1 & 1 & 1 & 1 & -1 & 1 & -1 & 1 & -1 & 1 & -1 &
1 & -1 & 1 \\
\chi_3 & 1 & -1 & 1 & -1 & 1 & 1 & 1 & -1 & 1 & -1 & 1 & 1 & -1 & 
1 & -1 & 1 & -1 \\
\chi_{13} & 2 & . & -2 & . & . & 2 & 2 & . & -2 & . & 2 & -2 & . & 2 & . & .
& . \\
\chi_{14} & 2 & . & -2 & . & . & 2 & 2 & . & -2 & . & 2 & 2 & . & -2 & . & .
& . \\
\chi_{17} & 8 & . & . & . & . & -8 & . & . & . & . & . & . & . & . & . & .
& . \\
\end{array} 
$$

\bigskip
{\bf Remark 2:} 
There is a real chance to solve the next regular case $d = 10$ in Sparla's
inequality. The question is whether there is a 2-Hamiltonian
4-manifold of genus 21 (i.e.  $\chi = 44$) in the 10-dimensional
cross polytope.
A 22-vertex triangulation of a manifold with 
exactly the same genus as a subcomplex
of the 11-dimensional cross polytope does exist.
If one could save two antipodal vertices by successive bistellar flips
one would have a solution.
The example with 22 vertices is 
defined by the orbits (of length 110 or 22, respectively) of the 4-simplices
$$
\langle 1\,3\,5\,7\,18 \rangle_{110}, 
\langle 1\,3\,5\,7\,21 \rangle_{110}, 
\langle 1\,3\,5\,8\,18 \rangle_{110}, 
\langle 1\,3\,5\,8\,21 \rangle_{110}, 
\langle 1\,3\,7\,18\,20 \rangle_{110},
\langle 1\,3\,6\,10\,15 \rangle_{22}$$
 under the permutation group of order 110 which is generated by
$$(1\ 16\ 7\ 22\ 13\ 5\ 19\ 12\ 3\ 18\ 10\ 2\ 15\ 8\ 21\ 14\ 6\ 20\ 11\ 4\ 17\ 9)$$ 
and 
$$(1\ 11\ 17\ 3\ 21)(2\ 12\ 18\ 4\ 22)(5\ 9\ 8\ 20\ 14)(6\ 10\ 7\ 19\ 13).$$
The central involution is 
$$(1\ 2)(3\ 4)(5\ 6)(7\ 8)(9\ 10)(11\ 12)(13\ 14)(15\ 16)(17\ 18)(19\ 20)(21\ 22)$$
which corresponds to the antipodal mapping in the cross polytope.
The $f$-vector of the example is $(22, 220, 1100, 1430, 572)$,
and the middle homology is 42-dimensional, the first and third
homology both vanish.
Hence it has ``genus'' 21 in the sense defined above.

\subsection*{5. Tightness and tautness}
The concept of tightness originates from differential geometry 
as the equality of
the (normalized) {\it total absolute curvature} of a submanifold with the
lower bound {\it sum of the Betti numbers} \cite{Kui}, \cite{Ba-Ku1}.
It is also a generalization of the concept of convexity since it roughly means
that an embedding of a submanifold is as convex as possible
according to its topology.
The usual definition is the following:
\subparagraph{Definition} (compare \cite{Kui})

An embedding $M \rightarrow \mathbb{E}^N$ of a compact manifold is called
{\sf tight}, if for any open or closed halfspace 
$E^N_+\subset \mathbb{E}^N$ the
induced homomorphism
$$H_*(M \cap E^N_+) \longrightarrow H_*(M)$$
is injective where $H_*$ denotes an appropriate homology theory
with coefficients in a certain field. 
The notion of {\it $k$-tightness} refers to the injectivity in the
low dimensions $H_i(M \cap E^N_+) \rightarrow H_i(M)$, $i = 0, \ldots, k$,
see \cite{Ku1}.
An equivalent formulation is that all nondegenerate height functions
are perfect functions, i.e., functions with a number of critical points
which equals the sum of the Betti numbers.
This definition applies to smooth and polyhedral embeddings.
A {\sf tight triangulation} is a triangulation of a manifold such that any
simplexwise linear embedding is tight \cite{Ku1}, \cite{K-Lu}.
Any $k$-Hamiltonian $2k$-manifold in the $d$-dimensional 
simplex is induced by a tight triangulation with $d+1$ vertices.
For a subcomplex of the boundary complex of a convex polytope the
tightness condition is often determined by purely combinatorial
conditions. In particular any $k$-Hamiltonian $2k$-manifold
in a $d$-polytope is tightly embedded into $d$-space \cite[4.1]{Ku1}. 
For any tight subcomplex $K$ of the boundary complex of a convex polytope
the following is a direct consequence of the definition above,
compare \cite[1.4]{Ku1}:
\subparagraph{Consequence}
{\it A facet of the polytope is either contained in $K$ or its 
intersection with $K$
represents a subset of $K$ (often called a topset)
which injects into $K$ at the homology level and which is again tightly
embedded into the ambient space.
In particular, any missing $(k+1)$-simplex in a $k$-Hamiltonian
subcomplex $K$ of a simplicial polytope represents a nonvanishing
element of the $k^{th}$ homology by the standard triangulation of the 
$k$-sphere.
}

\bigskip
For the similar notion of {\sf tautness} one has to replace halfspaces 
by balls (or ball complements) $B$
and height functions by distance functions, see \cite{CR}.
This applies only to smooth embeddings. In the polyhedral case it
has to be modified as follows:
\subparagraph{Definition} (suggested in \cite{Ba-Ku2})

A PL-embedding $M \rightarrow \mathbb{E}^N$ of a compact manifold 
with convex faces is called
{\sf PL-taut}, if for any open ball (or ball complement)
$B \subset \mathbb{E}^N$ the
induced homomorphism
$$H_*(M \cap {\rm span}(B_0)) \longrightarrow H_*(M)$$
is injective where $B_0$ denotes the set of vertices in  $M \cap B$,
and ${\rm span}(B_0)$ refers to the 
subcomplex in $M$ spanned by those vertices. 

\medskip
Obviously, any PL-taut embedding is also tight (consider very large balls),
and a tight PL-embedding is PL-taut provided that it is PL-spherical in 
the sense that all vertices are contained in a certain Euclidean
sphere. It follows that any tight and PL-spherical embedding is also PL-taut
\cite{Ba-Ku2}. 
\subparagraph{Corollary}
{\it Any tight subcomplex of a higher-dimensional regular 
simplex or cube or cross polytope is PL-taut.}

\medskip
In particular this implies that the class of PL-taut submanifolds
is much richer than the class of smooth taut submanifolds.
\subparagraph{Corollary}
{\it There is a tight and PL-taut simplicial embedding of
the connected sum of $7$ copies of $S^2 \times S^2$ into Euclidean $8$-space.
}

\bigskip
This follows directly from Theorem 2 by the embedding into
the 8-dimensional cross polytope.
In addition this example is centrally-symmetric.
There is a standard construction of tight embeddings of
connected sums of copies of $S^2 \times S^2$ but this
works in codimension 2 only,
polyhedrally as well as smoothly, see \cite[p.101]{Ba-Ku1}.
The cubical examples in \cite{KS} exist in arbitrary codimension
but they require a much larger ``genus'':
For a 2-Hamiltonian 4-manifold in the 8-dimensional cube one needs
an Euler characteristic $\chi = 64$ which corresponds
to a connected sum of 31 copies of $S^2 \times S^2$. 
The number of summands in this case grows exponentially with the
dimension of the cube.
For a 2-Hamiltonian 4-manifold in
the 8-dimensional simplex an Euler characteristic $\chi = 3$ is
sufficient. It is realized by the 9-vertex triangulation of ${\mathbb{C}}P^2$
\cite{Ku-Ba}, \cite{Ku-La1}.
One copy of $S^2 \times S^2$ cannot be a subcomplex of the 9-dimensional
simplex because such a 3-neighborly 10-vertex triangulation 
does not exist \cite{Ku-La1} even though it is one of the ``regular cases'' 
in the sense of the Heawood type integer condition in Section 1.
In general the idea behind is the following: A given $d$-dimensional
polytope requires a certain minimum ``genus'' of a 
$2k$-manifold to cover the full
$k$-dimensional skeleton of the polytope. For the standard polytopes
like simplex, $d$-cube and $d$-octahedron we have formulas for the ``genus''
which is to be expected but we don't yet have examples in all of the cases.

\medskip
The situation is similar with respect to the concept of tightness:
For any given dimension $d$ of an ambient space a certain ``genus''
of a manifold is required for admitting a tight and substantial
embedding into $d$-dimensional space. This is well understood in the
case of 2-dimensional surfaces \cite{Ku1}.
For ``most'' of the simply connected 4-manifolds a tight polyhedral 
embedding was constructed in \cite{Ku3}, without any especially intended
restriction concerning the essential codimension.
The optimal bounds in this case and in all the other higher-dimensional
cases still have to be investigated.

\subsection*{6. Centrally-symmetric triangulations of sphere products}
As far as the integer conditions of the ``regular cases''
are concerned, it seems to be plausible to ask for centrally-symmetric
triangulations of any sphere product $S^k \times S^l$
with a minimum number of $$n=2(k+l+2)$$
vertices. In this case each instance can be regarded
as a codimension-1-subcomplex of the boundary complex of the
$(k+l+2)$-dimensional cross polytope,
and that it can be expected to be $m$-Hamiltonian for $m = {\rm min}(k,l)$.
This is a kind of a simplicial Hopf decomposition of the $(k+l+1)$-sphere
by ``Clifford-tori'' of type $S^k \times S^l$.

\medskip
For $n \leq 20$ (i.e., for $k+l \leq 8$) a census of such triangulations with
a vertex-transitive automorphism group can be found in \cite{Lu},
compare \cite{K-Lu}.
Here all cases occur except for $S^4 \times S^2$ and $S^6 \times S^2$,
and all examples admit a dihedral group action of order $2n$.
So far an infinite series of examples
seems to be known only for $l=1$ and arbitrary $k$. This is the following:

\subparagraph{Proposition 4} (A centrally-symmetric and 1-Hamiltonian
$S^k \times S^1$  in $\partial \beta_{k+3}$)

{\it There is a centrally-symmetric triangulation of $S^k \times S^1$
with $n=2k+6$ vertices and with a dihedral automorphism group $D_n$.
Its induced embedding into the $(k+3)$-dimensional cross polytope
is tight and PL-taut.}

\bigskip
The construction is given in \cite{Ku-La2} with the notation
$M_k^{k+1}(n)$ (represented as the permcycle 
$[1^k2]$ there) as follows:
Regard the vertices as integers modulo $n$ and consider the
$\mathbb{Z}_n$-orbit of the $(k+2)$-simplex
$$\langle 0\,1\,2\, \cdots\, k\,(k+1)\,(k+2)\rangle.$$
This is a manifold with boundary (just an ordinary orientable 1-handle),
and its boundary is homeomorphic to $S^k\times S^1$.
All these simplices are facets of the cross polytope of dimension $k+3$
if we choose the labeling such that the diagonals are
$[x, x+k+3]$, $x \in \mathbb{Z}_n$. These diagonals do not
occur in the triangulation of the manifold, all other edges are contained.
Therefore we obtain a 1-Hamiltonian subcomplex of the $(k+3)$-dimensional
cross polytope. 
The central symmetry is the shift $x \mapsto x+k+3$ in $\mathbb{Z}_n$. The reflection $x \mapsto -x$ in $\mathbb{Z}_n$ is an extra automorphism.
In the case $k=1$ the group is even larger: It is of order 32.
This triangulated manifold is a hypersurface in $\partial \beta_{k+3}$,
it decomposes this $(k+2)$-sphere into two parts with
the same topology as suggested by the Hopf decomposition.

\medskip
The same generating simplex for the group $\mathbb{Z}_{m}$
with $m= 2k+5$ vertices leads to
the minimum vertex triangulation of $S^k \times S^1$ (for odd $k$)
or of the twisted product (for even $k$) which
is actually unique \cite{BD}, \cite{CS}.
For any $k\geq 2$ it realizes the minimum number of vertices for any manifold
of the same dimension which is not simply connected \cite{Br-Ku1}.
Other infinite series of triangulated sphere bundles over tori
are similarly given in \cite{Ku-La2}.

\medskip
It is not impossible that there will be direct generalizations
of Proposition 4 with infinite series of
analogous triangulations of $S^k\times S^3$, $S^k \times S^5$, \ldots,
at least for odd $k$,
and of $S^k \times S^k$, possibly for any $k$,
each with a dihedral and vertex-transitive group action and
Hamiltonian in the cross polytope. 
This is still work in progress.
Existence in the latter case of a $k$-Hamiltonian $S^k \times S^k$ 
with $n=4k+4$ vertices 
and $d=2k+2$ would give a positive answer to a conjecture by
F.~H.~Lutz \cite[p.85]{Lu}, and it 
would realize equality in Sparla's inequality in Section 3 for any $k$
since 
$$(-1)^k\frac{1}{2}\big(\chi - 2\big) = 1 = 
\frac{2k+1}{1} \cdot \frac{2k-1}{3} \cdot \frac{2k-3}{5} \cdot
\ \cdots \ \cdot \frac{1}{2k+1}.$$

\bigskip
Acknowledgement: We acknowledge support by the Deutsche Forschungsgemeinschaft
(DFG). This work was carried out as part of the DFG-project Ku 1203/5.

{\small

Institut f\"ur Geometrie und Topologie

Universit\"at Stuttgart

70550 Stuttgart

Germany

\medskip
effenberger@mathematik.uni-stuttgart.de

kuehnel@mathematik.uni-stuttgart.de
}

\begin{thebibliography}{45}
\bibitem{A} {\sc A.~Altshuler},
Manifolds in stacked $4$-polytopes, 
{\it J.\ Combin.\ Th.\ (A)} {\bf 10} (1971), 198--239

\bibitem{BD} {\sc B.~Bagchi} and {\sc B.~Datta},
Minimal triangulations of sphere bundles over the circle,\linebreak
{\it J.\ Combin. Th. (A)} {\bf 115} (2008), 737--752

\bibitem{BH} {\sc L.~W.~Beineke} and {\sc F.~Harary}, 
The genus of the $n$-cube, {\it Canad.\ J.\ Math.} {\bf 17} (1965), 494--496

\bibitem{Ba-Ku1} \textsc{T.~F.~Banchoff} and \textsc{W.~K\"uhnel}, 
Tight submanifolds, smooth and polyhedral, {\it Tight and Taut Submanifolds}
(T.~E.~Cecil and S.-S.~Chern, eds.), MSRI Publ.\ {\bf 32}, pp.\ 51--118,
Cambridge Univ.\ Press 1997 

\bibitem{Ba-Ku2} \textsc{T.~F.~Banchoff} and \textsc{W.~K\"uhnel}, 
Tight polyhedral models of isoparametric families, and PL-taut embeddings,
{\it Adv.\ Geom.} {\bf 7} (2007), 613--629

\bibitem{Bj-Lu} \textsc{A.~Bj{\"o}rner} and \textsc{F.~H.~Lutz}, 
Simplicial manifolds, bistellar flips and a 16-vertex triangulation of the {P}oincar\'e homology 3-sphere,
{\it Experiment.\ Math.} {\bf 9} (2000), 275--289

\bibitem{Br-Ku1} {\sc U.~Brehm} and {\sc W.~K\"uhnel}, 
Combinatorial manifolds with few vertices, 
{\it Topology} \textbf{26} (1987), 465--473

\bibitem{Br-Ku2} {\sc U.~Brehm} and {\sc W.~K\"uhnel}, 
$15$-vertex triangulations of an $8$-manifold, 
{\it Math.\ Ann.} \textbf{294} (1992), 167--193

\bibitem{CK} \textsc{M.~Casella} and \textsc{W.~K\"uhnel},
A triangulated $K3$ surface with the minimum  number of vertices,
{\it Topology} {\bf 40} (2001), 753--772

\bibitem{CR} \textsc{T.~E.~Cecil} and \textsc{P.~J.~Ryan},
{\it Tight and Taut Immersions of
Manifolds}, Pitman, Boston - London - Melbourne 1985

\bibitem{CS} \textsc{J.~Chestnut, J.~Sapir} and \textsc{E.~Swartz},
Enumerative properties of triangulations of sphere bundles over $S^1$,
{\it European J.\ Combin.} {\bf 29}1 (2008), 662--667

\bibitem{E} {\sc G.~Ewald}, Hamiltonian circuits in simplicial complexes,
{\it Geom.\ Ded.} {\bf 2} (1973), 115--125

\bibitem{GAP} \textsc{The GAP Group}, GAP - Groups, Algorithms and
  Programming, Version 4.2., Aachen, St.Andrews, 1999, see 
http://www.gap-system.org 

\bibitem{G-J} \textsc{E.~Gawrilow} and \textsc{M.~Joswig},
polymake: a framework for analyzing convex polytopes,
{\it Polytopes --- combinatorics and computation,} (Oberwolfach, 1997), 43--73, DMV Sem., {\bf 29}, Birkh{\"a}user, Basel, 2000.

\bibitem{GT} \textsc{R.~Z.~Goldstein} and \textsc{E.~C.~Turner},
A formula for Stiefel-Whitney homology classes,
{\it Proc.\ Amer.\ Math.\ Soc.} {\bf 58} (1976), 339--342

\bibitem{H} {\sc H.~Heesch}, {\it Gesammelte Abhandlungen},
Verlag Barbara Franzbecker 1986

\bibitem{J} \textsc{M.Joswig}, Computing invariants of simplicial manifolds, Preprint 2004, {\tt arXiv:math/0401176v1 [math.AT]}

\bibitem{JR} \textsc{M.~Jungerman} and \textsc{G.~Ringel},
The genus of the $n$-octahedron: Regular cases, 
{\it J.\ Graph Th.} {\bf 2} (1978), 69--75 

\bibitem{KT} \textsc{K.~M.~Koh} and \textsc{K.~L.~Teo},
The 2-Hamiltonian cubes of graphs, {\it J.\ Graph Th,} {\bf 13}
(1989), 737--747

\bibitem{Ku0} {\sc W.~K\"uhnel}, Manifolds in the skeletons of convex
polytopes, tightness, and generalized Heawood inequalities,
{\it POLYTOPES: Abstract, Convex and Computational,} 
Proc. Conf. Scarborough 1993 (T.~Bisztriczky et al., eds.),
241--247, NATO Adv. Study Inst. Ser. C, Math. Phys. Sci. {\bf 440},
Kluwer, Dordrecht 1994 

\bibitem{Ku1} {\sc W.~K\"uhnel}, 
{\it Tight Polyhedral Submanifolds and Tight Triangulations},
Lecture Notes in Mathematics \textbf{1612}, Springer, 
Berlin - Heidelberg - New York 1995

\bibitem{Ku2} {\sc W.~K\"uhnel}, Centrally-symmetric tight surfaces
and graph embeddings, {\it Beitr\"age Algebra Geom.} {\bf 37} (1996), 347--354

\bibitem{Ku3} {\sc W.~K\"uhnel}, Tight embeddings of simply connected
4-manifolds, {\it Doc.\ Math.} {\bf 9} (2004), 401--412

\bibitem{Ku-Ba} \textsc{W.~K\"uhnel} and \textsc{T.~F.~Banchoff}, 
The 9-vertex complex projective plane, 
{\it The Math.\ Intelligencer} \textbf{5}:3 (1983), 11--22

\bibitem{Ku-La1} \textsc{W.~K\"uhnel} and \textsc{G.~Lassmann}, 
The unique 3-neighborly
4-manifold with few vertices, {\it J.\ Combin.\ Th.\ (A)} \textbf{35} (1983), 
173--184

\bibitem{Ku-La2} \textsc{W.~K\"uhnel} and \textsc{G.~Lassmann}, 
Permuted difference cycles and triangulated sphere bundles,
{\it Discrete Math.} \textbf{162} (1996), 215--227

\bibitem{K-Lu} \textsc{W.~K\"uhnel} and \textsc{F.~H.~Lutz}, 
A census of tight triangulations, 
{\it Period.\ Math.\ Hung.} {\bf 39} (1999), 161--183

\bibitem{KS} \textsc{W.~K\"uhnel} and \textsc{C.~Schulz}, 
Submanifolds of the cube, {\it Applied Geometry and Discrete Mathematics,
The Victor Klee Festschrift} (P.Gritzmann and B.Sturmfels, eds.),
423--432, DIMACS Series in Discrete Math. and Theor.\ Comp.\ Sci. {\bf 4},
Amer.\ Math.\ Soc.\ 1991

\bibitem{Kui} \textsc{N.~H.~Kuiper}, 
Geometry in total absolute curvature theory, {\it Perspectives
in Mathematics}, Anniversary of Oberwolfach 1984 (W.\ J\"ager et al., eds.), 
377--392, Birkh\"auser, Basel - Boston - Stuttgart 1984

\bibitem{La-Sp} {\sc G.~Lassmann} and {\sc E.~Sparla}, A classification
of centrally-symmetric and cyclic $12$-vertex triangulations of 
$S^2 \times S^2$, {\it Discrete Math.} {\bf 223} (2000), 175--187

\bibitem{Lu}  \textsc{F.~H.~Lutz}, {\it Triangulated manifolds with few
vertices and vertex-transitive group actions}, Doctoral Thesis TU Berlin 1999, 
Shaker-Verlag, Aachen 1999

\bibitem{MMS} \textsc{P.~McMullen} and \textsc{E.~Schulte},
{\it Abstract Regular Polytopes}, Cambridge Univ. Press 2002 

\bibitem{N1} {\sc I.~Novik}, Upper bound theorems for homology manifolds,
{\it Israel J.\ Math.} {\bf 108} (1998), 45--82

\bibitem{N2} {\sc I.~Novik}, On face numbers of manifolds with symmetry,
{\it Adv.\ Math.} {\bf 192} (2005), 183--208

\bibitem{NS} {\sc I.~Novik} and {\sc E.~Swartz}, 
Socles of Buchsbaum modules, complexes and posets, Preprint 2007,
{\tt arXiv:0711.0783v1 [math.CO]}

\bibitem{R0} {\sc G.Ringel}, \"Uber drei Probleme am $n$-dimensionalen
W\"urfel und W\"urfelgitter, {\it Abh. Math. Sem. Univ. Hamburg} 
{\bf 20} (1955), 10--19 

\bibitem{R} {\sc G.~Ringel}, {\it Map Color Theorem}, Springer,
Berlin - Heidelberg - New York 1974 

\bibitem{Sa} {\sc N.~Saveliev}, {\it Lectures on the Topology of 
$3$-Manifolds. An Introduction to the Casson invariant}, 
de Gruyter, Berlin 1999 

\bibitem{Sch1} {\sc C.~Schulz}, {\it Mannigfaltigkeiten mit Zellzerlegung
im Randkomplex eines konvexen Polytops und verallgemeinerte Hamilton-Kreise},
Doctoral Dissertation Univ. Bochum 1974 

\bibitem{Sch2} {\sc C.~Schulz}, Polyhedral manifolds on polytopes,
Proc. Conf. Palermo 1993,
{\it Rend.\ Circ.\ Mat.\ Palermo (2) Suppl.} {\bf 35} (1993), 291--298

\bibitem{S} {\sc   C.~H.~S\'{e}quin}, 
Symmetrical Hamiltonian Manifolds on Regular 3D and 4D Polytopes, 
The Coxeter Day, Banff, Canada, Aug.\ 3, 2005, pp.\ 463--472, see\newline
http://www.cs.berkeley.edu/\~{}sequin/BIO/pubs.html  

\bibitem{Sp1} {\sc E.~Sparla}, An upper and a lower bound theorem for
combinatorial $4$-manifolds, {\it Discrete Comput.\ Geom.} {\bf 19} (1998), 
575--593 

\bibitem{Sp2} {\sc E.~Sparla}, A new lower bound theorem for
combinatorial $2k$-manifolds, {\it Graphs Combin.} {\bf 15} (1999), 
109--125 

\bibitem{Th} \textsc{G.~Thorbergsson}, Tight immersions of highly connected 
manifolds, {\it Comment.\ Math.\ Helv.} \textbf{61} (1986), 102--121

\end{thebibliography}
\end{document}